\documentclass[pdflatex,sn-mathphys-num]{sn-jnl}

\usepackage{graphicx}%
\usepackage{multirow}%
\usepackage{amsmath,amssymb,amsfonts}%
\usepackage{amsthm}%
\usepackage{mathrsfs}%
\usepackage[title]{appendix}%
\usepackage{xcolor}%
\usepackage{textcomp}%
\usepackage{manyfoot}%
\usepackage{booktabs}%
\usepackage{algorithm}%
\usepackage{algorithmicx}%
\usepackage{bm}
\usepackage{todonotes}
\usepackage{listings}%
\theoremstyle{thmstyleone}%
\theoremstyle{thmstyletwo}%
\theoremstyle{thmstylethree}%
\usepackage{float}
\usepackage{csquotes}
\usepackage[font={small,it}]{subcaption}
\usepackage{hyperref}
\usepackage[noend]{algpseudocode}
\algdef{SE}[SUBALG]{Indent}{EndIndent}{}{\algorithmicend\ }%
\algtext*{Indent}
\algtext*{EndIndent}
\newlength\myindent
\newcommand{\change}[1]{#1}

\begin{document}

\title{Augmented Lagrangian Solvers for Poroelasticity with Fracture Contact Mechanics}


\author*[1]{\fnm{Marius} \sur{Nevland}\href{https://orcid.org/0000-0002-1038-2052}{\textsuperscript{[ORCID]}}}\email{marius.nevland@uib.no}
\author[1]{\fnm{Inga} \sur{Berre}\href{https://orcid.org/0000-0002-0212-7959}{\textsuperscript{[ORCID]}}}\email{inga.berre@uib.no}
\author[1]{\fnm{Jakub} \sur{Wiktor Both}\href{https://orcid.org/0000-0001-6131-4976}{\textsuperscript{[ORCID]}}}\email{jakub.both@uib.no}
\author[1]{\fnm{Eirik} \sur{Keilegavlen}\href{https://orcid.org/0000-0002-0333-9507}{\textsuperscript{[ORCID]}}}\email{eirik.keilegavlen@uib.no}

\affil[1]{\orgdiv{Center for Modeling of Coupled Subsurface Dynamics, Department of Mathematics}, \orgname{University of Bergen}, \orgaddress{\street{Allégaten 41}, \city{Bergen}, \postcode{5020}, \country{Norway}}}

\abstract{In the subsurface, fractures and the surrounding porous rock can deform in interaction with fluid flow. Advanced mathematical models governing these coupled processes typically combine fluid flow, poroelasticity, and fracture contact mechanics. The resulting system of equations is complex and highly nonlinear. As a result, convergence issues with nonlinear solvers are common, \change{causing a bottleneck for the numerical solution of such models}.  

One source of difficulty for the nonlinear solvers comes from the fracture contact mechanics, due to its inherently nonsmooth character. \change{In addition, depending on the chosen constitutive model, the degree of nonlinearity is increased through coupling of flow and contact mechanics.} In this paper, we investigate solvers based on the augmented Lagrangian formulation of the frictional contact problem. This includes two classical solvers, namely the generalized Newton method (using complementarity functions) and the return map method \change{(equivalent to an Uzawa method). In addition, we propose a new solver that combines features of both approaches.}

Numerical experiments in two and three dimensions, \change{designed to simulate hydraulic stimulation of geothermal reservoirs}, are conducted to assess the performance of the solvers on problems of poromechanics with fracture contact mechanics. \change{The return map method has more difficulty handling the nonlinear coupling between flow and contact mechanics than the other solvers, in many cases not converging or using an excessive number of iterations. Our new combined solver performs the most robustly across the experiments, its performance being less sensitive to the value of the augmentation parameter than the other solvers.}}

\keywords{Fractured porous media, poroelasticity, contact mechanics, augmented Lagrangian, generalized Newton, return map, \change{Uzawa, projected Newton}}

\maketitle

\section{Introduction}

The mathematical modeling and simulation of coupled flow and deformation in fractured porous media has numerous important applications, such as in CO$_2$ storage and geothermal energy production. The tightly coupled physical processes, combined with the discontinuous nature of fractures and the structural complexity of fracture networks, result in significant mathematical challenges, from a modeling as well as a numerical perspective.

An overview of different conceptual approaches to model fractured porous media can be found in \cite{fracture_review,viswanathan2022fluid,vaezi2025review}. In this work, we specifically consider discrete fracture-matrix models, in which the fractures are explicitly represented. Specifically, the matrix (i.e. the porous medium surrounding the fractures) is modeled as linearly poroelastic \cite{coussy} with flow governed by Darcy's law. For the fracture, Darcy flow is also considered, while deformation is modeled by contact mechanics including Coulomb friction \cite{kikuchi_oden}. There have been several recent works on such models \cite{berge,Stefansson_2021,stefansson2023flexible,NOVIKOV2022111598,garipov2016,GARIPOV2019104075,BEAUDE2023116124,droniou_masson,SALIMZADEH2018212}, focusing on the development of discretization methods adapted to the discrete fracture-matrix model. In this context, for the mechanical subproblem, a common approach is to use finite element methods \cite{garipov2016,GARIPOV2019104075,BEAUDE2023116124,droniou_masson,SALIMZADEH2018212}, but, motivated by the coupling with flow, finite volume methods are also increasingly being applied \cite{berge,Stefansson_2021, stefansson2023flexible,NOVIKOV2022111598}. The flow equations are \change{commonly} discretized by finite volume methods \cite{berge,Stefansson_2021,stefansson2023flexible,NOVIKOV2022111598,garipov2016,GARIPOV2019104075,BEAUDE2023116124,droniou_masson}, but finite element methods have also been used \cite{SALIMZADEH2018212}. In addition, considerable efforts have been made to stabilize the mechanical subproblem \cite{FRANCESCHINI2020113161,BEAUDE2023116124,uzawa_bubble}.

To apply these discretization schemes to complex setups that reflect the physical states and fractured structures encountered in field applications, 
there is a need for nonlinear solvers for discrete fracture-matrix models for poromechanics, and scalable linear solvers.
Linear solvers, based on block preconditioning that is tailored to both the individual physical processes and their couplings, 
have been developed for both finite element~\cite{franceschini2019block,franceschini2022scalable,franceschini2022scalable2}
and finite volume~\cite{zabegaev2025efficient,zabegaev2025block} discretizations.
In terms of nonlinear solvers, a challenging aspect is the inherent nonlinear and nonsmooth character of frictional contact problems~\cite{kikuchi_oden,WriggersPeter2006CCM}. Several classes of nonlinear solvers that were originally developed for frictional contact mechanics have been extended to account for poromechanics. These solvers are built around different ways of formulating the contact conditions, and includes solvers based on the penalty \cite{NOVIKOV2022111598, garipov2016}, augmented Lagrangian \cite{berge,Stefansson_2021,BEAUDE2023116124,SALIMZADEH2018212} and Nitsche formulations \cite{GARIPOV2019104075, BEAUDE2023116124} of the frictional contact problem. However, the poromechanical coupling introduces new challenges for the nonlinear solvers. First, the extra couplings make for larger and more complex algebraic systems.
Second, depending on model choices, the systems might contain other nonlinearities that are triggered by the contact mechanics, such as shear dilation of fractures and fracture permeability modeled by a cubic law \cite{Stefansson_2021,stefansson2023flexible}.
As a result, the design of robust nonlinear solvers for these problems remains a challenge that must be addressed to enable reliable simulations with advanced physical models on complex fracture geometries.

\change{Nonlinear solvers based on the augmented Lagrangian formulation of frictional contact mechanics typically follows one of two distinctly different approaches.} The first approach is based on a reformulation of the variational inequalities of contact mechanics as the zero sets of certain non-smooth functions called complementarity functions. This enables the use of a generalized Newton method with appropriately defined generalized Jacobians to handle the system’s non-smoothness. This approach was first introduced by Alart and Curnier \cite{ALART1991353}, and several different ways of formulating the complementarity functions have since been derived \cite{klarbring,hueber}. The second approach, introduced by Simo and Laursen \cite{SIMO199297}, uses a fixed-point iteration on the augmented Lagrange multipliers, resulting in \change{an implicit} return map method that can also be interpreted as an \change{implicit} Uzawa algorithm \cite{RenardSurvey}. \change{Both approaches are well-established in the contact mechanics literature, and a systematic comparison of these two solvers (in addition to other nonlinear solvers) has been conducted in \cite{RenardSurvey} for problems of a linearly elastic body in contact with a rigid obstacle, but not considering fractured media.}

\change{The generalized Newton method based on complementarity functions has been used in several of the aforementioned works to solve problems of coupled flow and deformation including fracture contact mechanics \cite{berge, Stefansson_2021,BEAUDE2023116124}. The implicit return map method, using the augmented Lagrangian formulation, appears to have been less used to solve such problems. It was used together with a pure penalty formulation in \cite{garipov2016,NOVIKOV2022111598}. In \cite{nejati2016finite}, an augmented Lagrangian method was developed for purely mechanical contact problems in fractured media, where a return map was used to augment the fracture gap, rather than the Lagrange multipliers. The result is an algorithm that is distinctly different from the one by Simo and Laursen \cite{SIMO199297}. It was subsequently used in \cite{SALIMZADEH2018212} to solve a coupled thermo-hydro-mechanical model.}

\change{In the present work, we perform a systematic study and comparison of three different nonlinear solvers based on the augmented Lagrangian formulation. These include the generalized Newton method of Alart and Curnier \cite{ALART1991353}, using the complementarity functions of Hüeber et al. \cite{hueber}, and the implicit return map method of Simo and Laursen \cite{SIMO199297}. A systematic comparison of these two solvers, in the vein of \cite{RenardSurvey}, is lacking for problems of coupled flow and deformation including fracture contact mechanics. The third solver in our study is new, and it combines features of both the previous two solvers. Specifically, it is a quasi-Newton method that incorporates the return map as a postprocessing step after every generalized Newton iteration. Conceptually, this solver can be interpreted as a type of projected Newton method, which is a method used to solve constrained optimization problems \cite{projected_newton}, where each Newton iteration is followed by a return map (projection) onto the feasible set. A comparable approach has been employed by Jha et al. \cite{jha2014coupled} to address purely mechanical subproblems within a fixed-stress solver for coupled multiphase poromechanics in fractured media. However, to the authors' knowledge, the present work is the first time this type of solver has been applied to solve coupled flow and mechanics in a monolithic manner.}

\change{Moreover, we show a mathematical connection between the generalized Newton and the implicit return map methods. Specifically, we will
demonstrate that the implicit return map method is equivalent to solving a sequence of nonlinear systems with
regularized complementarity functions, where each of these systems can be solved by the generalized
Newton method. To the authors' knowledge, the literature lacks a detailed derivation of this connection, which is key to understanding how our new, combined solver differs from both the generalized Newton and the implicit return map methods.}

To illustrate the convergence challenges and evaluate the performance of the generalized Newton method, the \change{implicit} return map method, and our new nonlinear solver, we carry out a series of numerical experiments in both two and three dimensions. These experiments are inspired by hydraulic stimulation of geothermal reservoirs, in which critically stressed fractures are stimulated by fluid injection, resulting in fracture slip and corresponding shear dilation. \change{In the simulations, we use different constitutive model choices for the relationship between slip, aperture dilation and fracture permeability, to gradually increase the degree of nonlinear couplings between flow and mechanics, highlighting how such couplings can contribute to convergence difficulties. We also investigate a range of values for the so-called augmentation parameter; this is a numerical parameter inherent to the augmented Lagrangian formulation, which must be chosen by the user. This is relevant to study, as the value of this parameter can affect the convergence behavior of the nonlinear solvers, as shown in the study done in \cite{RenardSurvey} on purely mechanical contact problems.}

The remainder of the paper is structured as follows. In Section \ref{governing_eqn}, we present the mathematical model for coupled flow and deformation in fractured porous media. Section \ref{sec:nonlinear_solvers} introduces the generalized Newton method and the \change{implicit} return map method, adapted to the governing equations, and establishes a formal connection between them. \change{Our new nonlinear solver, which combines features of both previous solvers, is subsequently introduced.} In Section \ref{num_exp}, we evaluate the performance of all three nonlinear solvers through a series of numerical experiments. Concluding remarks are provided in Section \ref{sec:conclusion}.

\section{Mathematical model}

\label{governing_eqn}

We consider a discrete fracture-matrix model, where the fractures are modeled as co-dimension one objects, leading to a mixed-dimensional poromechanical model including fracture contact mechanics.  The presentation of the model in this section closely follows that of \cite{stefansson2023flexible}; see also \cite{berge,porepy,Stefansson_2021}. We first give an overview of the mixed-dimensional geometry in section \ref{sec:mixed_geometry}. In section \ref{matrix_eqns}, conservation laws for mass and momentum balance, with corresponding constitutive laws, are presented. The equations of fracture contact mechanics are presented in section \ref{fracture_contact}. Finally, in section \ref{initial_bc} we define the boundary conditions and provide a summary of the primary variables and equations.

\subsection{Mixed-dimensional geometry} \label{sec:mixed_geometry}

The fractured porous medium is described as a collection of subdomains $\Omega_i$ of varying dimension $d_i \in \{0,...,D\}$, with $D \in \{2,3\}$. The matrix subdomain will have dimension $D$, while fractures are modeled as co-dimension one objects, and will accordingly have dimension $D-1$. Moreover, intersections of fractures will have dimension $D-2$, and in the case of $D=3$, one also considers intersections of fracture intersection lines as zero-dimensional subdomains. The width of a fracture is characterized by its \change{volumetric} aperture $a_i$. This is further generalized by introducing the specific volume,
\begin{equation} \label{eq:specific_volume}
   \change{\mathcal{V}_i=a_i^{D-d_i},}
\end{equation}
accounting for dimension reduction for a subdomain $\Omega_i$ with any dimension $d_i \leq D$.

\begin{figure}
[h]\centering
    {\scalebox{0.3}{\includegraphics{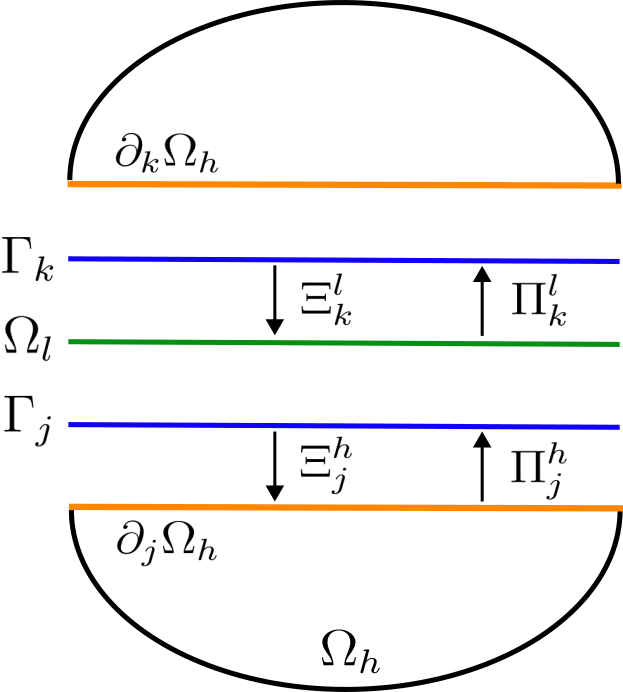}}}
    \caption{Illustration of the mixed-dimensional geometry. The higher-dimensional subdomain $\Omega_h$ is coupled to the lower-dimensional subdomain $\Omega_l$ through the interfaces $\Gamma_j$ and $\Gamma_k$.  Note that $\partial_j \Omega_h, \partial_k \Omega_h, \Gamma_j, \Gamma_k$ and $\Omega_l$ all coincide geometrically. Figure adapted from \cite{Stefansson_2021}}
    \label{hamburger}
\end{figure}

A pair of subdomains one dimension apart will always be connected by one or more interfaces $\Gamma_j, \Gamma_k$. This is done to facilitate the coupling of processes between subdomains of varying dimension. For such a pair of subdomains, it will be convenient to denote the higher-dimensional subdomain by $\Omega_h$ and the lower-dimensional one by $\Omega_l$. As illustrated in figure \ref{hamburger}, the interfaces $\Gamma_j, \Gamma_k$ will geometrically coincide with $\Omega_l$ and the internal boundaries $\partial_j \Omega_h, \partial_k \Omega_h$ of $\Omega_h$. It is sometimes necessary to project relevant quantities from a subdomain to an interface or vice versa. The projection from a subdomain $\Omega_i$ to an interface $\Gamma_j$ is denoted by $\Pi_j^i$, while the reverse operation is denoted by $\Xi_j^i$.

Subscripts will be used to denote on which subdomain or interface a quantity is defined. A generic subdomain and interface are denoted by subscripts $i$ and $j$, respectively. If two interfaces are involved, we will use subscript $k$ for the second interface. For certain equations involving couplings between a higher- and lower-dimensional subdomain, we use subscript $h$ for the higher dimension and $l$ for the lower dimension. 

\subsection{Conservation laws and constitutive relations} \label{matrix_eqns}

Mixed-dimensional conservation laws establish the foundation for the communication  across the different dimensions.
The fluid mass conservation equation for a general subdomain $\Omega_i$ of dimension $d_i \in \{ 0,...,D \}$ reads:
\begin{equation} \label{eq:fluid_mass_balance}
    \frac{\partial}{\partial t}(\mathcal{V}_i \rho_i^f \phi_i) + \nabla \cdot (\mathcal{V}_i \rho_i^f \bm{v}_i) - \sum_{j \in \hat{S}_i}\Xi_j^i (\mathcal{V}_j \rho_j^f v_j) = \psi_i, 
\end{equation}
where $\rho_i^f, \rho_j^f$ are the subdomain and interface fluid densities, respectively, $\phi_i$ is the porosity, $\bm{v}_i$ and $v_j$ are the subdomain and interface volumetric fluid fluxes, $\psi_i$ is a source term, $\mathcal{V}_j := \Pi_j^i\mathcal{V}_i$ is the interface specific volume, and the set $\hat{S}_i$ contains all interfaces connecting to higher-dimensional neighbors of $\Omega_i$. The second term on the left hand side of~\eqref{eq:fluid_mass_balance} vanishes for $d_i=0$, as zero-dimensional domains do not have mass fluxes, while the third term on the left hand side of~\eqref{eq:fluid_mass_balance} vanishes for $d_i=D$, in which case there are no higher-dimensional neighboring subdomains.

Assuming a quasi-static process, the momentum balance equation of the matrix subdomain ($d_i=D$) reads:
\begin{equation} \label{mom_bal}
    -\nabla \cdot \bm{\sigma} = \bm{F},
\end{equation}
with $\bm{\sigma}$ being the total poroelastic stress tensor and $\bm{F}$ denoting body forces. 

Next, we define constitutive relations for the various quantities defined above, most of which are derived from \cite{coussy}. The volumetric fluid flux is modeled by Darcy's law:
\begin{equation} \label{darcy}
    \bm{v}_i = -\frac{\bm{\mathcal{K}}_i}{\eta}(\nabla p_i-\rho_i^f \bm{g}),
\end{equation}
where $\eta$ is the viscosity (assumed to be constant), $p_i$ is the pressure, $\bm{\mathcal{K}}_i$ is the permeability tensor and $\bm{g}$ the gravitational acceleration vector. In order to define the constitutive laws for the permeabilities, we first define $\mathbf{I}_D$ and $\mathbf{I}_{D-1}$ to be the identity matrices of dimension $D$ and $D-1$, respectively. The permeability tensor $\bm{\mathcal{K}}_i$ is assumed to be spatially homogeneous in the matrix, hence it is a constant multiple $k$ of the identity matrix:
\begin{equation} \label{eq:matrix_permeability}
   \bm{\mathcal{K}}_i = k\mathbf{I}_D \ \ , \ \ d_i=D, 
\end{equation}
while in fracture subdomains, it is given by a cubic law \cite{permeability}:
\begin{equation} \label{eq:tang_perm}
\bm{\mathcal{K}}_i = \frac{\change{\mathcal{A}_i^2}}{12}\mathbf{I}_{D-1} \ \ , \ \ d_i=D-1,
\end{equation}
\change{where $\mathcal{A}_i$ denotes the hydraulic aperture, which we distinguish from the volumetric aperture $a_i$. Constitutive laws for both apertures will be defined in section \ref{fracture_contact}, after some new notation has been introduced.} Finally, intersection permeability is computed as the average permeability of the intersecting fractures. Solid density is assumed to be constant, while the density of the slightly compressible fluid is given by (we omit the subscript $i$ for readability):
\begin{equation} \label{density_cons}
    \rho^f = \rho_\mathrm{ref}^f \exp{(\gamma(p-p_\mathrm{ref}))},
\end{equation}
where $\gamma$ is the compressibility of the fluid. The subscript ``ref'' is used to denote reference values of variables. The interface flux $v_j$ is given by an averaged Darcy's law \cite{martin}:
\begin{equation} \label{dim_red_darcy}
    v_j = -\frac{\mathcal{K}_j}{\eta}\left(\frac{2}{\Pi^l_j a_l}(\Pi_j^l p_l-\Pi_j^h p_h)- \rho_j^f \bm{g} \cdot 
 \Pi_j^h \bm{n}_h \right).
\end{equation}
Here, $\bm{n}_h$ is the outward normal vector of the higher-dimensional neighboring subdomain. The interface permeability $\mathcal{K}_j$ is also modeled by a cubic law, inheriting the \change{hydraulic} aperture from its lower-dimensional neighboring subdomain,
\begin{equation} \label{interface_permeability}
   \mathcal{K}_j=\frac{\Pi_j^l\change{\mathcal{A}_l^2}}{12},
\end{equation}
while the interface fluid density $\rho_j^f$ is defined by an upstream weighting based on the direction of $v_j$:
\begin{equation} \label{upstream}
    \rho_j^f = \begin{cases}
        \Pi^h_j\rho_h^f \ \ \ &\text{if} \ \ v_j>0, \\
        \Pi^l_j \rho_l^f \ \ \ \ &\text{if}  \ \ v_j \leq 0.
    \end{cases}
\end{equation}
The total poroelastic stress tensor is given by an extended Hooke's law that also accounts for fluid pressure:
\begin{equation} \label{stress_tensor}
    \bm{\sigma} = G(\nabla \bm{u} + \nabla \bm{u}^T) + \lambda_{\text{Lamé}}\text{tr}(\nabla \bm{u})\mathbf{I}_D - \alpha p\mathbf{I}_D.
\end{equation}
Here, $\bm{u}$ denotes the displacement, $G$ is the shear modulus, $\lambda_{\text{Lamé}}$ is Lamé's first parameter, $\alpha$ is the Biot coefficient and $\text{tr}(\cdot)$ denotes the trace of a matrix. The porosity of fractures and intersections is set to one, while the matrix porosity depends on displacement and pressure as follows:
\begin{equation} \label{porosity_cons}
    \phi = \phi_{\mathrm{ref}} + \frac{(\alpha-\phi_{\mathrm{ref}})(1-\alpha)}{\lambda_{\text{Lamé}}+\frac{2}{3}G}(p-p_{\mathrm{ref}}) + \alpha(\nabla \cdot \bm{u}).
\end{equation}
The body force in the momentum balance equation \eqref{mom_bal} is given as the combined gravitational forces of the solid and the fluid:
\begin{equation} \label{body_force}
    \bm{F} = (\phi \rho^f +(1-\phi)\rho^s)\bm{g},
\end{equation}
where $\rho^s$ is the solid density.
\subsection{Fracture contact mechanics}
\label{fracture_contact}
The equations of contact mechanics are imposed on pairs of $D$-dimensional matrix subdomains $\Omega_h$ and $D-1$-dimensional fracture subdomains $\Omega_l$ that are connected through two interfaces $\Gamma_j, \Gamma_k$, as illustrated in figure \ref{hamburger}. Some new variables are introduced \change{at this point, namely} the fracture contact traction $\bm{\lambda}_l$, defined on the fracture domain $\Omega_l$, and the displacements $\bm{u}_j, \bm{u}_k$ of the interfaces $\Gamma_j, \Gamma_k$; \change{these are treated as additional primary unknowns to be solved for.} For readability purposes, we will omit the subscript $l$ from the contact traction. The forces on the two fracture surfaces must be balanced according to Newton's third law. Taking into account that the pressure from the fluid in the fracture contributes to the total fracture traction, the following balance of forces are imposed at the interfaces:
\begin{equation} \label{force_interface}
    \begin{aligned}
        \Pi_j^l (\bm{\lambda} - p_l \mathbf{I}_D \cdot \bm{n}_l) &= \Pi_j^h \sigma_h \cdot \bm{n}_h \ \ \ \ \ \ \text{on} \ \Gamma_j, \\
         \Pi_k^l (\bm{\lambda} - p_l \mathbf{I}_D \cdot \bm{n}_l) &= -\Pi_k^h \sigma_h \cdot \bm{n}_h \ \ \ \ \text{on} \ \Gamma_k, \\
    \end{aligned}
\end{equation}
where $\bm{n}_l$ denotes the normal vector of the fracture, which is defined as $\bm{n}_l=\Xi_j^l \Pi_j^h \bm{n}_h$. In other words, it is chosen to equal $\bm{n}_h$ on the $j$ side. A vector $\bm{v}$ defined on a fracture may be decomposed into its normal and tangential components with respect to the fracture as follows:
\begin{equation}
    v_n=\bm{v} \cdot \bm{n}_l, \ \ \ \bm{v}_{\tau}=\bm{v}-v_n \bm{n}_l.
\end{equation}
The jump in interface displacements across $\Omega_l$ is defined as $\relax[[\bm{u}]]=\Xi_k^l \bm{u}_k-\Xi_j^l \bm{u}_j$. The relative movement of the fracture is restricted by a nonpenetration condition, which can be written as:
\begin{equation} \label{nonpen}
\begin{aligned}
    \relax[[\bm{u}]]_n-g &\geq 0, \\
    \lambda_n &\leq 0, \\
    \lambda_n([[\bm{u}]]_n-g) &= 0,
\end{aligned}
\end{equation}
where $g$ is the gap function, which should have the property that $\relax[[\bm{u}]]_n-g=0$ when the surfaces are in contact, and $\relax[[\bm{u}]]_n-g > 0$, otherwise. The constitutive law for $g$ depends on modeling choices. In this work we use a formula that incorporates shear dilation of fractures. Due to small scale surface roughness of a fracture, displacement along the fracture can cause shear dilation, resulting in the fracture being mechanically closed but hydraulically open. This effect is incorporated into the model by defining a dilation angle $\psi$, and choosing the gap function to be 
\begin{equation} \label{eq:gap}
    g(\bm{u})=\tan(\psi)\lVert [[\bm{u}]]_{\tau} \rVert.
\end{equation}
This definition results in $g$ having the aforementioned properties. Next, we require the tangential surface traction to satisfy the Coulomb friction law in the case of contact, while it should vanish in the case of no contact:
\begin{equation} \label{friction_law}
\begin{cases}
\text{If} \ \lambda_n=0 \implies \bm{\lambda}_{\tau}=0. \\
\text{Else}: \\
    \ \ \begin{aligned}
        \lVert \bm{\lambda}_{\tau} \rVert &\leq -\change{\mu}\lambda_n, \\
         \lVert \bm{\lambda}_{\tau} \rVert &< -\change{\mu}\lambda_n \implies [[\dot{\bm{u}}]]_{\tau}=\bm0, \\
         \lVert \bm{\lambda}_{\tau} \rVert &= -\change{\mu}\lambda_n \implies \exists \ \zeta \in \mathbb{R}^+ : [[\dot{\bm{u}}]]_{\tau} = \zeta \bm{\lambda}_{\tau} ,
    \end{aligned}
\end{cases}
\end{equation}
where $\change{\mu}$ is the coefficient of friction and $[[\dot{\bm{u}}]]_{\tau}$ is the time derivative of the tangential displacement jump. We will for convenience refer to the entirety of equation \eqref{friction_law} as the Coulomb friction law. Finally, we define the constitutive \change{laws for the volumetric and hydraulic apertures. In fracture subdomains, these apertures are set to equal the normal component of the interface displacement jump, with additional constant reference apertures $a_{\mathrm{ref}}, \mathcal{A}_{\mathrm{ref}}$ taking into account small-scale roughnesses in the undeformed state}:
\begin{equation} \label{eq:vol_aperture}
    a_{\change{i}} = a_{\mathrm{ref}} + [[\bm{u}]]_n \change{\ \ , \ \ d_i=D-1,}
\end{equation}
\begin{equation} \label{eq:hyd_aperture}
    \change{\mathcal{A}_{i} = \mathcal{A}_{\mathrm{ref}} + [[\bm{u}]]_n \ , \ \ d_i=D-1.}
\end{equation}
For intersection subdomains, \change{the apertures are computed as the average apertures of the intersecting fractures}:
\begin{equation} \label{intersection_vol_aperture}
    a_i = \frac{1}{\lvert \hat{S}_i \rvert}\sum_{j \in \hat{S}_i}\Xi_j^i \Pi_j^h a_h \change{\ \ , \ \ d_i<D-1.}
\end{equation}
\begin{equation} \label{intersection_hyd_aperture}
    \change{\mathcal{A}_i = \frac{1}{\lvert \hat{S}_i \rvert}\sum_{j \in \hat{S}_i}\Xi_j^i \Pi_j^h \mathcal{A}_h \ \ , \ \ d_i<D-1.}
\end{equation}
\change{The volumetric aperture $a_i$ enters the fracture flow model through the specific volume \eqref{eq:specific_volume}, while the hydraulic aperture $\mathcal{A}_i$ governs the fracture permeability through the cubic law \eqref{eq:tang_perm}. Both cases result in a strong nonlinear coupling. We hypothesize a strong impact of the degree of coupling and nonlinearity on the performance of the nonlinear solvers. This will be investigated in the numerical experiments of section \ref{num_exp}, by a systematic comparison of different linearizations of the apertures, where either one or both of the volumetric and hydraulic apertures are kept fixed to reduce the degree of nonlinearity.}
\subsection{Boundary conditions and summary} \label{initial_bc}

Boundary conditions are imposed on both internal and external boundaries. On the internal boundaries $\partial_j \Omega_i$, we require continuity of normal mass fluxes and displacements (the latter only applies for $d_i=D$):
\begin{align}
\mathcal{V}_i\rho_i^f \bm{v}_i \cdot \bm{n}_i &= \Xi_j^i \mathcal{V}_j \rho_j^f v_j, \\
\bm{u}_i&=\Xi_j^i \bm{u}_j.
\end{align}
On external boundaries, we impose either Dirichlet or Neumann boundary conditions. The mechanical boundary data is accordingly either prescribed values of displacement $\bm{u}$ or traction $\sigma \cdot \bm{n}$, and fluid boundary data is either prescribed values of pressure $p$ or mass flux $\rho \bm{v} \cdot \bm{n}$. Finally, we impose zero fluid mass flux on immersed fracture tips, supported by the discussion in~\cite{both2024}.

A summary of the system of equations and primary variables to be solved for is outlined in table \ref{eqns_vars}.

\begin{table}[ht!]
\centering
\begin{tabular}{ | l | l | l | }
\hline
Equations & Variables & Subdomains or interfaces \\
\hline
Mass balance, \eqref{eq:fluid_mass_balance} & Pressure $p_i$ & $\Omega_i$ for all $d_i$ \\
Momentum balance, \eqref{mom_bal} & Displacement $\bm{u}_i$ & $\Omega_i$ for $d_i=D$ \\
Interface Darcy flux, \eqref{dim_red_darcy} & Interface Darcy flux $v_j$ & $\Gamma_j$ for all $d_j$ \\
Interface force balance, \eqref{force_interface} & Interface displacement $\bm{u}_j$ & $\Gamma_j$ for $d_j=D-1$ \\
Nonpenetration condition, \eqref{nonpen} & Contact traction $\bm{\lambda}$ & $\Omega_i$ for $d_i=D-1$ \\
Coulomb friction law \eqref{friction_law} & & $\Omega_i$ for $d_i=D-1$ \\
\hline
\end{tabular}
\caption{Summary of equations and primary variables, as well as the subdomains or interfaces on which they are defined. The equation and variable on a specific row are both defined on the same subdomains or interfaces, which is specified in the same row.}
\label{eqns_vars}
\end{table}

\section{Augmented Lagrangian solvers} \label{sec:nonlinear_solvers}

In this section, we will present three different nonlinear solvers for the system of equations presented in the previous section. All three solvers are based on the augmented Lagrangian formulation of the frictional contact problem. Two of these solvers are classical, namely the generalized Newton method and the implicit return map, while the third solver is new. To ease the presentation and stress their agnostic character to spatial discretizations, the solvers will be presented in a semi-discrete setting, where the equations are assumed to be discretized in time, using the implicit Euler method, but continuous in space. We will give some details in section \ref{sec:alg_full_disc} of the particular spatial discretizations used in the numerical experiments of the present work.

\textbf{Notation}. \change{The full set of primary unknowns is collected into a vector $\bm{x}$, consisting of the pressure, displacement, interface Darcy flux, interface displacement and contact traction, as outlined in table \ref{eqns_vars}. Nonlinear iteration indices will be denoted by a superscript, i.e. $\bm{x}^{(k)}$. Moreover, we will separate the contact conditions, \eqref{nonpen} and \eqref{friction_law}, from the rest of the equation system. This is done because the implicit return map method enforces the contact conditions differently from the other two solvers, while the remaining equations are treated the same by all solvers. We denote by $G$ the residual system in which the contact conditions have been excluded. Hence, the equation}
\begin{equation} \label{residual_system}
    G(\bm{x})=\bm{0},
\end{equation}
\change{constitutes the semi-discrete formulations of equations \eqref{eq:fluid_mass_balance}, \eqref{mom_bal}, \eqref{dim_red_darcy} and \eqref{force_interface}, or in other words, every equation of table \ref{eqns_vars} with the exception of the nonpenetration conditon \eqref{nonpen} and Coulomb friction law \eqref{friction_law}.}

\subsection{The generalized Newton method} \label{sec:generalized_newton}
The nonpenetration condition \eqref{nonpen} and Coulomb friction law \eqref{friction_law} can be reformulated as the zero sets of certain nonsmooth functions called complementarity functions. Such a reformulation can be derived from a theoretical consideration of the augmented Lagrangian formulation of the contact problem, as shown in \cite{ALART1991353}.
There are several ways of formulating these complementarity functions \cite{ALART1991353,klarbring,hueber}. In the present work, we employ the formulation of Hüeber et al. \cite{hueber}, in which it is shown that, for every augmentation parameter $c>0$, the nonpenetration condition (\ref{nonpen}) is equivalent to:
\begin{equation} \label{compl_func_normal}
    \mathcal{C}_n :=\lambda_n+\max(0,-\lambda_n-c([[\bm{u}]]_n-g(\bm{u}))) = 0.
\end{equation}
Similarly, for every $c>0$, the Coulomb friction law (\ref{friction_law}) is equivalent to:
\begin{equation} \label{compl_func_tangential}   \mathcal{C}_{\tau} :=\chi \bm{\lambda}_{\tau}+(1-\chi)\left[b(\bm{\lambda}_{\tau}+c[[\dot{\bm{u}}]]_{\tau})-\max(b,\lVert \bm{\lambda}_{\tau}+c[[\dot{\bm{u}}]]_{\tau} \rVert)\bm{\lambda}_{\tau}\right] = \bm{0}.
\end{equation}
\change{Here $b=-\mu\lambda_n$ is the friction bound and $\chi$ is a characteristic function that equals 1 if $\lvert b\rvert \leq 0$ and 0 if $\lvert b\rvert > 0$, thus enforcing $\bm{\lambda}_{\tau}=\bm{0}$ in the case of no contact. Note that, as the equations in this section are assumed to be discretized in time, the slip velocity $[[\dot{\bm{u}}]]_{\tau}$ now refers to the time increment of the tangential displacement jump. In other words, we have $[[\dot{\bm{u}}]]_{\tau}=[[\bm{u}]]_{\tau}-[[\bm{u}]]_{\tau}^{\text{prev}}$, where ``prev'' denotes the previous time step, and the time step size is implicitly subsumed in the augmentation parameter $c$}.

By replacing \eqref{nonpen} and \eqref{friction_law} with \eqref{compl_func_normal} and \eqref{compl_func_tangential}, the entire equation system of table \ref{eqns_vars} may be written compactly in residual form. Using the previously defined residual system~\eqref{residual_system}, the full system of equations to be solved at a specific time step reads 
\begin{equation} \label{algebraic_system}
    F(\bm{x}):= \begin{bmatrix} G(\bm{x}) \\ \mathcal{C}_n(\bm{x}) \\ \mathcal{C}_\tau(\bm{x}) \end{bmatrix} = \bm{0}.
\end{equation}
This system is only differentiable almost everywhere (in the sense of the Lebesgue measure), due to the non-smoothness of the complementarity functions. However, the non-differentiable points can easily be handled by employing generalized derivatives, using the definition found in Clarke \cite{clarke}. A set of generalized Jacobians for the system \eqref{algebraic_system} at a point $\bm{x}^*$ is defined as follows:
\begin{equation}
    \partial F(\bm{x}^*) = \text{conv}\left\{ \lim_{\bm{x}_i \to \bm{x}^*, \  \bm{x}_i \in \Omega_F} DF(\bm{x}_i) \right\},
\end{equation}
where $\Omega_F$ is the set of points on which $F$ is differentiable, with corresponding Jacobian $DF$, and ``conv'' denotes the convex hull. Using these generalized Jacobians, a corresponding generalized Newton method can be defined, \change{with the resulting algorithm being equivalent to a primal-dual active set method \cite{hintermuller2002primal,hueber}. We describe this generalized Newton method below, and abbreviate it by GNM.}
\begin{algorithm}[H]
\caption{Generalized Newton Method (GNM)}
\label{alg:generalized_newton}
\begin{algorithmic} 
\State 1. Pick an initial guess $\bm{x}^{(0)}$. \change{Choose tolerances TOL1 and TOL2.} Set $k=0$.

\State 2. Solve for $\bm{x}^{(k+1)}$ the Newton linearization of system (\ref{algebraic_system}) at $\bm{x}^{(k)}$:

\begin{equation} \label{eq:newton_iteration}
    V^{(k)}(\bm{x}^{(k+1)}-\bm{x}^{(k)}) = -F(\bm{x}^{(k)}),
\end{equation}

\State \hspace{0.36cm} where $V^{(k)} \in \partial F(\bm{x}^{(k)})$.

\State 3. Convergence check:
\Indent
\If{\change{$\lVert\bm{x}^{(k+1)}-\bm{
x}^{(k)}\rVert<\text{TOL1}$ and $\lVert F(\bm{x}^{(k+1)}\rVert<\text{TOL2}$}}
    \State Exit algorithm.
\Else
    \State Set $k=k+1$ and go to 2.
\EndIf
\EndIndent

\end{algorithmic}
\end{algorithm}

Let us briefly comment on the augmentation parameter $c$ appearing in the complementarity functions. Since the two complementarity functions are an exact reformulation of the nonpenetration condition \eqref{nonpen} and Coulomb friction law \eqref{friction_law} for all $c>0$, this parameter may \change{in theory} be chosen freely without affecting the accuracy of the solution. \change{However, on a practical level,} the choice of the parameter can greatly affect the performance of the generalized Newton method (this is also the case for the other nonlinear solvers considered in this work), \change{and this will be an important point in the numerical experiments of section \ref{num_exp}}.

\change{Note also that the augmentation parameter can be chosen differently for the normal \eqref{compl_func_normal} and tangential \eqref{compl_func_tangential} complementarity functions. However, we will for simplicity assume throughout this paper that it is equal for the two complementarity functions; this will also hold for the subsequent algorithms of sections \ref{sec:classic_return} and \ref{sec:newton_return_map}. In our numerical experiments, we did not find any significant enough advantages in using two different parameters to warrant the added complexity and bookkeeping.}

\subsection{The implicit return map method} \label{sec:classic_return}

Return map methods are commonly used in computational plasticity \cite{comp_plasticity, computational_inelasticity}, and have also been applied to frictional contact problems \cite{GIANNAKOPOULOS1989157,wriggers_return_map,SIMO199297}. The defining feature of such methods is that for every iteration, a trial value is proposed for a variable subject to constraints, and violation of these constraints is checked for this value. If the trial value lies outside the feasible set, it is projected back to the boundary of the set; this projection is what is called the return map. This methodology was first applied to the augmented Lagrangian formulation of frictional contact mechanics in \cite{SIMO199297}. 

The central idea of the augmented Lagrangian-based return map is to define a trial value of the contact traction at every iteration, where the trial value is defined as the sum of the traction at the previous iteration and a penalty term monitoring the constraints on the normal and tangential displacement jumps. The penalty term, together with the return map, should augment the contact traction until the constraints on the displacement jumps converge, in which case the traction should also converge. In view of this, the trial values should be, at some iteration $k$:
\begin{subequations}
\label{eq:trial_tractions_uzawa}
\begin{align}
\lambda_n^{\text{trial}} &= \lambda_n^{(k)} + \change{c^{(k)}}([[\bm{u}^{(k+1)}]]_n-g(\bm{u}^{(k+1)})), \\
  \bm{\lambda}_{\tau}^{\text{trial}} &= \bm{\lambda}_{\tau}^{(k)}+\change{c^{(k)}}[[\dot{\bm{u}}^{(k+1)}]]_{\tau},    
\end{align} 
\end{subequations}
for some \change{$c^{(k)}>0$, which may be chosen differently for each iteration}. These formulas result in the contact traction converging in the case when $[[\bm{u}]]_n-g \to 0$ and $[[\dot{\bm{u}}]]_{\tau} \to \bm{0}$ (stick), while the remaining cases (open, slip) are handled by the return map. Some new notation is introduced at this point: We separate the contact traction $\bm{\lambda}$ from the rest of the unknowns, and denote the \change{remaining} unknowns by $\bm{r}$. As we have previously defined the total collection of unknowns as $\bm{x}$, we may write $\bm{x}=(\bm{r}, \bm{\lambda})$. The residual system \eqref{residual_system} is accordingly rewritten as $G(\bm{r,\bm{\lambda}})=\bm{0}$. Using this notation, we formulate the return map method below, which is \change{equivalent to} the method from \cite{SIMO199297}. We \change{refer to} it as the \textit{implicit return map method} (IRM).

\begin{algorithm}[H] 
\caption{Implicit Return Map Method (IRM)}
\label{alg:simo_laursen}
\begin{algorithmic}
\State 1. Pick an initial guess $\bm{x}^0=(\bm{r}^{0}, \bm{\lambda}^{0})$. \change{Choose tolerances TOL1 and TOL2.} Set $k=0.$
\State 2. Solve the following system for $\change{\bm{x}^{(k+1)}}=(\bm{r}^{(k+1)},\bm{\lambda}^{(k+1)}):$
\begin{equation} \label{eq:residual_simo}
G(\bm{r}^{(k+1)},\bm{\lambda}^{(k+1)})=\bm{0},
\end{equation}
\Indent
\State where $\bm{\lambda}^{(k+1)}$ is given by the following return map:
\begin{subequations} \label{eq:simo_return_map}
    \begin{align}
        \lambda_n^{(k+1)}&=\begin{cases}
            \lambda_n^{\text{trial}} \qquad &\text{if} \ \ \lambda_n^{\text{trial}}\leq 0, \\
            0  &\text{if} \ \ \lambda_n^{\text{trial}}> 0,
        \end{cases} \label{eq:simo_return_norm} \\
        \bm{\lambda}_{\tau}^{(k+1)} &= \begin{cases}
            \bm{\lambda}_{\tau}^{\text{trial}} \qquad\qquad &\text{if} \ \ \lVert\bm{\lambda}_{\tau}^{\text{trial}}\rVert\leq b^{(k+1)}, \\
            b^{(k+1)}\frac{\bm{\lambda}_{\tau}^{\text{trial}}}{\lVert \bm{\lambda}_{\tau}^{\text{trial}}\rVert} &\text{if} \ \ \lVert\bm{\lambda}_{\tau}^{\text{trial}}\rVert> b^{(k+1)},
        \end{cases} \label{eq:simo_return_tang}
    \end{align}
\end{subequations}
\State where $b^{(k+1)}=-\change{\mu}\lambda_n^{(k+1)}$, \change{and $\lambda_n^{\text{trial}}$, $\bm{\lambda}_{\tau}^{\text{trial}}$ are given by \eqref{eq:trial_tractions_uzawa}.}
\EndIndent
\State 3. Convergence check:
\Indent
\If{\change{$\lVert \bm{x}^{(k+1)}-\bm{x}^{(k)}\rVert<\text{TOL1}$ and $\lVert F(\bm{x}^{(k+1)})\rVert<\text{TOL2}$}}
    \State Exit algorithm.
\Else
    \State Set $k=k+1$ and go to 2.
\EndIf
\EndIndent
\end{algorithmic}
\end{algorithm}

\change{As alluded to by the naming of the algorithm, the return map is only performed implicitly, due to the presence of $\bm{\lambda}^{(k+1)}$ in \eqref{eq:residual_simo}. This means that \eqref{eq:residual_simo} and \eqref{eq:simo_return_map} must be solved simultaneously. We note also that algorithm \ref{alg:simo_laursen} defines a nested loop; each outer iteration $k$ consists of solving the system defined by \eqref{eq:residual_simo} and \eqref{eq:simo_return_map}. This is a nonlinear system, which must be solved in a separate inner loop. The outer loop can be interpreted as the augmentation procedure of the Lagrange multipliers (contact tractions), in addition to re-solving the residual system with the updated Lagrange multipliers. Note also that the convergence check is done on the system \eqref{algebraic_system} defined earlier, thus ensuring that the contact constraints are adequately satisfied by checking the residuals of the complementarity functions \eqref{compl_func_normal} and \eqref{compl_func_tangential}.}

We will now provide a connection between the return maps \eqref{eq:simo_return_norm} and \eqref{eq:simo_return_tang} and the complementarity functions \eqref{compl_func_normal} and \eqref{compl_func_tangential}, which leads to a practical way of solving the system \eqref{eq:residual_simo} and \eqref{eq:simo_return_map}. The key observation is that the return maps can be reformulated using max-functions. Additionally, we need to define the characteristic function $\chi^{(k+1)}$, equal to one if $b^{(k+1)}=0$ and zero if $b^{(k+1)}>0$. Note that $b^{(k+1)}$ will always be non-negative, since the return map ensures that $\lambda_n^{(k+1)}\leq0$. Using this notation, equations \eqref{eq:simo_return_norm} and \eqref{eq:simo_return_tang} can be rewritten as:
\begin{subequations} \label{eq:return_map_fpi}
  \begin{align}
    \lambda_n^{(k+1)} &= -\max(0,-\lambda_n^{\text{trial}}), \\
    \bm{\lambda}_{\tau}^{(k+1)} &= \frac{(\chi^{(k+1)}-1)b^{(k+1)}\bm{\lambda}_{\tau}^{\text{trial}}}{(\chi^{(k+1)}-1)\max(b^{(k+1)},\lVert \bm{\lambda}_{\tau}^{\text{trial}}\rVert)+\chi^{(k+1)}}.
\end{align}  
\end{subequations}
The characteristic function is needed to avoid division by zero in the case when $b^{(k+1)}=\lVert\bm{\lambda}_{\tau}^{\text{trial}}\rVert=0$. If we now further rewrite \eqref{eq:return_map_fpi} to residual form and substitute the definitions of the trial tractions, we obtain:
\begin{subequations} \label{eq:regularized_comp}
    \begin{align}
        &\lambda_n^{(k+1)} +\max(0, -\lambda_n^{(k)} - \change{c^{(k)}}([[\bm{u}^{(k+1)}]]_n-g(\bm{u}^{(k+1)}))) = 0, \label{eq:reg_norm} \\
        &\chi^{(k+1)} \bm{\lambda}_{\tau}^{(k+1)}+(1-\chi^{(k+1)})\left[b^{(k+1)}(\bm{\lambda}_{\tau}^{(k)}+\change{c^{(k)}}[[\dot{\bm{u}}^{(k+1)}]]_{\tau})-\max(b^{(k+1)},\lVert \bm{\lambda}_{\tau}^{(k)}+\change{c^{(k)}}[[\dot{\bm{u}}^{(k+1)}]]_{\tau}\rVert)\bm{\lambda}_{\tau}^{(k+1)}\right] = \bm{0}. \label{eq:reg_tang}
    \end{align}
\end{subequations}
Equations \eqref{eq:reg_norm} and \eqref{eq:reg_tang} are nearly identical to the complementarity functions \eqref{compl_func_normal} and \eqref{compl_func_tangential}; the only differences are the terms involving the traction at the previous iteration, $\bm{\lambda}^{(k)}$, which is a known quantity in these equations. One can interpret \eqref{eq:reg_norm} and \eqref{eq:reg_tang} as regularized contact conditions, in the sense that they convert the contact conditions from graphs to functions, due to the $\bm{\lambda}^{(k)}$-terms. Hence, IRM is equivalent to solving a sequence of systems with regularized contact conditions. More precisely, for a fixed $\bm{\lambda}^{(k)}$ one may define the regularized complementarity functions
\begin{align*}
\mathcal{C}_{n}^{(k)} &:=
\lambda_n +\max(0, -\lambda_n^{(k)} - \change{c^{(k)}}([[\bm{u}]]_n-g(\bm{u}))), \label{eq:reg_norm2}
\\
\mathcal{C}_{\tau}^{(k)} &:= \chi\bm{\lambda}_{\tau}+(1-\chi)\left[b(\bm{\lambda}_{\tau}^{(k)}+\change{c^{(k)}}[[\dot{\bm{u}}]]_{\tau})-\max(b,\lVert \bm{\lambda}_{\tau}^{(k)}+\change{c^{(k)}}[[\dot{\bm{u}}]]_{\tau}\rVert)\bm{\lambda}_{\tau}\right].
\end{align*}
Iteration $k$ of IRM then amounts to solving for $\bm{x}^{(k+1)}$ the following regularized version of system \eqref{algebraic_system}:
\begin{equation} \label{eq:regularized_system}
    \begin{bmatrix}        G(\bm{x}^{(k+1)}) \\
        \mathcal{C}_{n}^{(k)}(\bm{x}^{(k+1)}) \\
        \mathcal{C}_\tau^{(k)}(\bm{x}^{(k+1)})
    \end{bmatrix} = \bm{0}.
\end{equation}
These regularized systems may be solved by GNM (algorithm \ref{alg:generalized_newton}), and this is how IRM has been implemented in practice in the numerical experiments of section \ref{num_exp}. We note also that the return map \eqref{eq:simo_return_map} is equivalent to updating the contact tractions (or augmented Lagrange multipliers) by fixed-point iteration, according to \eqref{eq:return_map_fpi}. This corresponds to the classical \textit{Uzawa} algorithm \cite{WriggersPeter2006CCM}. Hence, IRM can also be interpreted as an Uzawa algorithm. Again we emphasize that the fixed-point iteration is only performed implicitly in IRM, resulting in an implicit Uzawa algorithm. It is known that for purely mechanical contact problems, implicit Uzawa algorithms converge for a larger range of augmentation parameters than their explicit counterparts \cite{Kunisch_Stadler_2005}, at the cost of needing to solve a larger and more complicated system.

The augmentation parameter $c$ takes on a different meaning for IRM, compared to GNM. Here it directly controls the strength of the regularization, and as $c \to \infty$, the regularized complementarity functions \eqref{eq:reg_norm} and \eqref{eq:reg_tang} converge to the exact functions \eqref{compl_func_normal} and \eqref{compl_func_tangential}. \change{However, the augmentation procedure on the contact tractions ensures that IRM will in theory converge to the non-regularized solution for any $c>0$, just like GNM. That being said, IRM still differs from GNM in that} the convergence speed is highly dependent on the value of $c$, with lower values requiring more iterations of the outer loop to converge. This is a well-known feature of Uzawa algorithms \cite{RenardSurvey,Kunisch_Stadler_2005}. \change{Speedup of the outer loop may be achieved by gradually increasing $c$ with increasing outer loop iterations. Since $c$ need not be taken to infinity to obtain convergence, one may define an upper limit for its value and stop the increase once this limit is reached.}

\subsection{A generalized Newton method with a return map} \label{sec:newton_return_map}

\change{In this section, we will present a new augmented Lagrangian solver, which combines features of both previous solvers. The central idea of this solver is to perform the same generalized Newton iterations as GNM, but after each iteration, a return map of the contact tractions is performed before proceeding to the next iteration. This is motivated by the fact that Newton solvers are well-known to often display divergent behavior (e.g. cycling) outside the radius of quadratic convergence. Performing a return map after each iteration might help to stabilize the solver in such cases. The solver is described below, using the notation introduced in section \ref{sec:classic_return}. We refer to it as the \textit{generalized Newton method with a return map}, and abbreviate it by GNM-RM.}

\begin{algorithm}[H]
\caption{Generalized Newton Method with a Return Map (GNM-RM)}
\label{newton_return_map}
\begin{algorithmic}
\State 1. Pick an initial guess $\bm{x}^0=(\bm{r}^0, \bm{\lambda}^0)$. \change{Choose tolerances TOL1 and TOL2.} Set $k=0$.

\State 2. \change{Perform one generalized Newton iteration \eqref{eq:newton_iteration} on the system \eqref{algebraic_system} to get \\ \hspace{0.4cm}$\bm{x}^{\text{trial}}=(\bm{r}^{(k+1)}, \bm{\lambda}^{\text{trial}}).$}

\State 3. Convergence check:
\Indent
\If{\change{$\lVert \bm{x}^{\text{trial}}-\bm{x}^{(k)}\rVert<\text{TOL1}$ and $\lVert F(\bm{x}^{\text{trial}})\rVert<\text{TOL2}$}}
    \State Exit algorithm.
\Else
    \State \textbf{Return map in normal direction:}
    \State $\lambda_n^{(k+1)}=\begin{cases}
            \lambda_n^{\text{trial}} \ \ \ \ \text{if} \ \ \lambda_n^{\text{trial}}\leq 0, \\
            0 \ \ \ \ \ \ \ \ \  \text{if} \ \ \lambda_n^{\text{trial}}> 0.
        \end{cases}$           
    \State \textbf{Return map in tangential direction:}
    \State Define the new friction bound $b^{(k+1)}=-\change{\mu}\lambda_n^{(k+1)}$.
    \vspace{0.1cm}
    \State $\bm{\lambda}_{\tau}^{(k+1)}= \begin{cases}       \bm{\lambda}_{\tau}^{\text{trial}} \ \ \ \ \ \ \ \ \ \ \ \ \ \ \ \ \  \text{if} \ \ \lVert\bm{\lambda}_{\tau}^{\text{trial}}\rVert\leq b^{(k+1)}, \\
    b^{(k+1)}\frac{\bm{\lambda}_{\tau}^{\text{trial}}}{\lVert \bm{\lambda}_{\tau}^{\text{trial}}\rVert} \ \ \ \ \ \ \  \text{if} \ \ \lVert\bm{\lambda}_{\tau}^{\text{trial}}\rVert> b^{(k+1)}.
    \end{cases}$
    \State \change{Update guess to $\bm{x}^{(k+1)}=(\bm{r}^{(k+1)}, \bm{\lambda}^{(k+1)})$, set $k=k+1$ and go to 2.}
\EndIf
\EndIndent
\end{algorithmic}
\end{algorithm}

\change{We note that the return map in GNM-RM is explicitly computed after each generalized Newton iteration, which is in contrast to IRM, where the return map is implicitly defined through the regularized complementarity functions. GNM-RM is hence distinctly different from an Uzawa method like IRM, and should instead be interpreted as a quasi-Newton method.}

\change{GNM-RM is conceptually similar to the so-called projected Newton method \cite{projected_newton} from the optimization literature. This method solves optimization problems with inequality constraints by following each Newton iteration with a return map (projection) onto the feasible set. However, the return map plays a conceptually different role in GNM-RM, compared to both the projected Newton method of \cite{projected_newton}, as well as IRM. For the latter two methods, the return map (projection) is performed in order to ensure the fulfillment of the inequality constraints. For GNM-RM, on the other hand, the fulfillment of the inequality constraints, i.e the contact conditions, is already handled by performing generalized Newton iterations on the complementarity functions \eqref{compl_func_normal} and \eqref{compl_func_tangential}. Hence, the return map of GNM-RM functions merely as stabilization of GNM.}

\subsection{Finite volume discretization} \label{sec:alg_full_disc}

We conclude this section with a brief summary of the spatial discretizations used in this work, and the corresponding practical implementation of algorithms \ref{alg:generalized_newton}-\ref{newton_return_map}. We follow the discretization schemes described in \cite{stefansson2023flexible} and the references therein. \change{In short, we discretize the conservation laws by a family of cell-centered finite volume schemes. For the fluid mass balance equation, the standard multi-point flux approximation is used, while for the momentum balance, we use the multi-point stress approximation, first introduced in \cite{nordbotten2014cell}. Advective terms are discretized by an upwinding scheme.} The grids are simplicial and conforming to the fractures, in the sense that lower-dimensional cells coincide with higher-dimensional faces. 

In this fully discrete formulation, the residual system \eqref{residual_system} becomes a nonlinear system of algebraic equations, while the complementarity functions are enforced cell-wise, resulting in another system of nonlinear algebraic equations. The implementation of the generalized Newton method on the full system (or on the regularized systems in the case of IRM) is now straightforward. In GNM-RM, the return maps are performed cell-wise.

\section{Numerical experiments} \label{num_exp}

In this section, we will compare the performance of the three nonlinear solvers presented in section \ref{sec:nonlinear_solvers}, GNM, IRM and GNM-RM, through a series of two- and three-dimensional numerical simulations. Specifically, we simulate pressurized fluid injection in a fractured reservoir, mimicking hydraulic stimulation of a fractured geothermal reservoir in crystalline rock. \change{Before moving on to the experiments themselves, we will first outline the main objectives of these experiments, and provide some technical details of the implementation.}

\subsection{Main objectives} \label{sec:main_objectives}

\change{A central focus of the experiments is to test the impact of gradually increasing the strength of the nonlinear couplings between flow and contact mechanics on the nonlinear solver performance. Since the three nonlinear solvers under consideration are all built around the contact problem, and since GNM and IRM are already well-established for solving purely mechanical contact problems \cite{RenardSurvey}, it is important to investigate to which degree the nonlinear coupling of the contact problem to flow contributes to convergence issues. We identify the nonlinear dependence of the volumetric and hydraulic apertures on the normal displacement jump, $a_i=a_{\text{ref}}+[[\bm{u}]]_n$ and $\mathcal{A}_i=\mathcal{A}_{\text{ref}}+[[\bm{u}]]_n$, as central. Hence, we will consider three different models, with varying degrees of linearizations of these aperture laws. We will first consider a model where both the effective fracture volume and fracture permeabilities are considered constant, using $a_i=a_{\text{ref}}$ and $\mathcal{A}_i=\mathcal{A}_{\mathrm{ref}}$. Note that this only (partially) linearizes the flow equations, while the mechanical gap function \eqref{eq:gap} will still be active. Next, we consider an intermediate problem with volumetric changes, but a constant fracture permeability. This means that a constant hydraulic aperture $\mathcal{A}_i=\mathcal{A}_{\text{ref}}$ is used, but the volumetric aperture formula still depends on the normal displacement jump, $a_i=a_{\text{ref}}+[[\bm{u}]]_n$. Finally, we include also the cubic law on the fracture permeability and consider the fully coupled model, using the displacement-dependent formulas for both the volumetric and hydraulic apertures. These three models are summarized in table \ref{tab:models_experiments}.}

\begin{table}[ht!]
\centering
\begin{tabular}{ | l | l | l | }
\hline
Model & Feature & Aperture formulas used \\
\hline
A & Constant fracture volume and permeability & $a_i=a_{\text{ref}}, \  \mathcal{A}_i=\mathcal{A}_{\text{ref}}$ \\
B & Constant fracture permeability & $a_i=a_{\text{ref}}+[[\bm{u}]]_n$, \ $\mathcal{A}_i=\mathcal{A}_{\text{ref}}$\\
C & Full model & $a_i=a_{\text{ref}}+[[\bm{u}]]_n$, \ $\mathcal{A}_i=\mathcal{A}_{\text{ref}}+[[\bm{u}]]_n$ \\
\hline
\end{tabular}
\caption{The three models, with varying degrees of nonlinear couplings of flow and mechanics, used in the simulations.}
\label{tab:models_experiments}
\end{table}

\change{Another goal of the experiments is to test the sensitivity of the nonlinear solvers to the value of the augmentation parameter $c$. This is relevant to study, as the value of this parameter can greatly influence the convergence behavior of the solvers. This is demonstrated in \cite{RenardSurvey} for purely mechanical contact problems, where GNM and IRM were shown to display sensitive behavior with respect to the value of $c$. In particular, for some of the examples, the set of $c$-values for which the solver converged was quite erratic and unpredictable.} 

\change{We will consider $c$-values in a range centered around the value of 1 GPa/m. We choose $c$-values around this order, as we are aiming to balance the traction- and displacement-terms of the complementarity functions, and it is empirically found that the ratio $\bm{\lambda}/[[\bm{u}]]$ is typically in the range from $10^{-2}$ to 1 GPa/m in these simulations. This balancing of terms has been recommended in several articles on purely mechanical contact problems \cite{hueber,gitterle}. For IRM, the value of $c$ will be gradually increased during each iteration of the outer loop; this is done because the convergence of IRM will be slower for smaller $c$-values, as mentioned in section \ref{sec:classic_return}. Specifically, if $c$ denotes the starting value, our update strategy is $c^{(k)}=\max(10^3, 10^k \cdot c)$, with $k$ denoting the outer loop iteration of IRM. The value is reset back to $c$ for each new time step. This feature of smaller $c$-values being slower to converge is unique to Uzawa methods like IRM, and hence, it does not apply to GNM and GNM-RM, and we simply set $c^{(k)}=c$ in these two algorithms.}

\change{Our final objective, which is related to the previous two objectives, is to assess and compare the overall performance of the three nonlinear solvers. As a related sub-objective, we will also report and discuss the performance of the linear solver, details of which can be found in the next section. The performance of each nonlinear solver will be assessed by reporting the total number of nonlinear iterations needed for convergence. By this we mean the total number of linear systems solved, summed over all time steps. For the convergence criteria, we use the same absolute tolerance of $10^{-8}$ for both the nonlinear iteration increments and residuals, both of which are measured in a Euclidean norm that is weighted by the cell volumes. We note that for IRM, the tolerance of $10^{-8}$ is used for both the outer loop and for solving the regularized system in the inner loop. An alternative would be to use a looser tolerance in the inner loop, and adaptively tighten the tolerance with increasing outer loop iterations.  However, in our experiments, we did not experience any significant advantages by using a looser tolerance in the inner loop.}

\subsection{Technical details of the implementation} \label{sec:technical_details}

All of the simulations are run using \change{version 1.11 of} the mixed-dimensional simulation framework PorePy \cite{porepy}. \change{The grids are generated by Gmsh \cite{gmsh} (version 4.13.1), while linear solver routines are imported from PETSc \cite{balay2023petsc,DALCIN20111124} using \textit{petsc4py} (version 3.24.0)}. Runscripts for reproducing the results of this section can be found in the Docker image in \cite{zenodo}. The precise vertex coordinates of the fractures used in the simulations are reported in the supplementary material of this article.

\change{The linear systems are solved by GMRES, using absolute and relative tolerances of $10^{-11}$, together with the preconditioner developed in \cite{zabegaev2025efficient}, which is tailored to mixed-dimensional poroelasticity with contact mechanics. Specifically, it is a Schur complement based preconditioner using the algebraic multigrid method, with the ILU($n$) smoother being used for the mechanics subproblem. The value of $n$ has been slightly varied across the experiments, to ensure the most robust performance of the linear solver. Specifically, we have used $n=2$ for the two-dimensional simulations and $n=0$ for the three-dimensional ones. In some cases, however, the linear solver failed to converge; in these cases, the adaptive time-stepping scheme described below is employed.}

\change{Some measures are taken to improve the conditioning of the nonlinear systems. First, to avoid having to deal with negative apertures in the non-converged state, we use a cutoff for the volumetric and hydraulic apertures that do not allow them to be less than their respective residual apertures. Thus, we are in practice using the formulas $a_i=\max(a_{\text{ref}},a_{\text{ref}}+[[\bm{u}]]_n)$ and $\mathcal{A}_i=\max(\mathcal{A}_{\text{ref}},\mathcal{A}_{\text{ref}}+[[\bm{u}]]_n)$. Next, we have observed empirically that the return map of GNM-RM will quite frequently result in cases where $\lVert \bm{\lambda}_{\tau}^{(k+1)}+c[[\dot{\bm{u}}^{(k+1)}]]_{\tau}\rVert=b^{(k+1)}$, i.e. the two arguments of the maximum-function of \eqref{compl_func_tangential} are equal. In this case, GNM-RM was found to be more robust if the Jacobian corresponding to the second argument is chosen, and so this is done in all simulations. For GNM and IRM, neither of the two options showed a significant increase or decrease in the robustness of the solvers. Therefore, we arbitrarily choose the Jacobian corresponding to the second argument in case of equal arguments also for GNM and IRM. Finally, we use an appropriate scaling of the units. It is found empirically that setting the unit of mass to be $10^9$ kg results in the primary variables being close to the same order of magnitude, and generally improves the robustness of all solvers. Hence, for all simulations, the input data is scaled by this unit of mass.} 

\change{An adaptive time-stepping scheme is used in the simulations, inspired by the one used by Simunek et al. \cite{simunek2005hydrus}. The time-stepping scheme has two features. First, if the nonlinear solver converged, the next time step size is adapted based on the number of nonlinear iterations needed for convergence. If the number of iterations is within a specified range, the step size stays the same; if it is above this range, the step size is decreased by a factor of 0.7, while if it is below this range, the step size is increased by a factor of 3. The range will slightly vary in the examples. For the two-dimensional simulations of section \ref{sec:two_dim_sim}, we set the range to $(4,20)$, not including the endpoints. For the three-dimensional simulations of section \ref{subsec:three_dim_ex}, we increase the lower endpoint of the range to avoid excessive computational time, setting it to $(8,20)$. In section \ref{sec:three_dim_difficult} we run even more challenging three-dimensional simulations, in which case we shrink the range even further, to $(10,20)$.}

\change{The second feature of the time-stepping scheme is that, if the nonlinear solver fails to converge, a recomputation of the solution is attempted with a smaller step size. We use three metrics to assess non-convergent behavior: Setting the maximal number of allowed nonlinear iterations per nonlinear system to be 30, setting the maximal value of the absolute residual norm to be $10^5$, and setting the maximal number of allowed GMRES iterations to be 300. If any of these cases occur, the step size is decreased by a factor of one half, and a recomputation of the nonlinear system is attempted with this updated step size. If IRM fails to converge in the inner loop, the whole outer loop is restarted with the smaller time step. We allow a step size no smaller than a tenth of the initial time step, and if the simulation still fails to converge for this step size, it is terminated. We also terminate the simulation after 6 consecutive failed recomputation attempts.}

\change{Line search methods are often used as globalization strategies for nonlinear solvers. We experimented with combining the nonlinear solvers with an Armijo line search, but found that this did not improve the robustness of the solvers for the problem setups considered herein. A similar behavior has been reported for certain problems of purely mechanical contact mechanics \cite{Acary2018}. Hence, we have chosen to run all of the simulations without a line search strategy. Finally, the material constants used in the simulations are listed in table \ref{material_constants}, in unscaled units.}

\begin{table}[ht]
\centering
\begin{tabular}{ l c c c }
\hline
 Parameter & Symbol & Value & Units \\ 
\hline
 Biot coefficient & $\alpha$ & 0.8 & - \\  
 Friction coefficient & $\change{\mu}$ & 0.5 & - \\
 Matrix permeability & $k$ & $1.0 \cdot 10^{-15}$ & m$^2$ \\
 Reference porosity & $\phi_{\mathrm{ref}}$ & $1.0 \cdot 10^{-2}$ & - \\
 Shear modulus & $G$ & $1.7 \cdot 10^{10}$ & Pa \\
Lamé's first parameter & $\lambda_{\text{Lamé}}$ & $1.111 \cdot 10^{10}$ & Pa \\
 Compressibility & $\gamma$ & $0.4 \cdot 10^{-9}$ & Pa$^{-1}$ \\
 Reference fluid density & $\rho^f_{\mathrm{ref}}$ & $1.0 \cdot 10^3$ & kgm$^{-3}$ \\
 Solid density & $\rho^s$ & $2.7 \cdot 10^3$ & kgm$^{-3}$ \\
 Viscosity & $\eta$ & $1.0 \cdot 10^{-3}$ & Pa \ s \\
 Residual \change{volumetric} aperture & $a_{\mathrm{ref}}$ & $5.0 \cdot 10^{-4}$ & m \\
 \change{Residual hydraulic aperture} & \change{$\mathcal{A}_{\mathrm{ref}}$} & \change{$5.0 \cdot 10^{-4}$} & \change{m} \\
 Dilation angle & $\psi$ & \change{5} & degrees\\
 Reference pressure & $p_{\mathrm{ref}}$ & $2.0 \cdot 10^7 $ & Pa
\end{tabular}
\caption{Material parameters used in the simulations.}
\label{material_constants}
\end{table}

\subsection{Two-dimensional simulations of fluid injection} \label{sec:two_dim_sim}

\begin{figure}
[H]\centering
{\scalebox{1.0}
{\includegraphics[width=\textwidth]{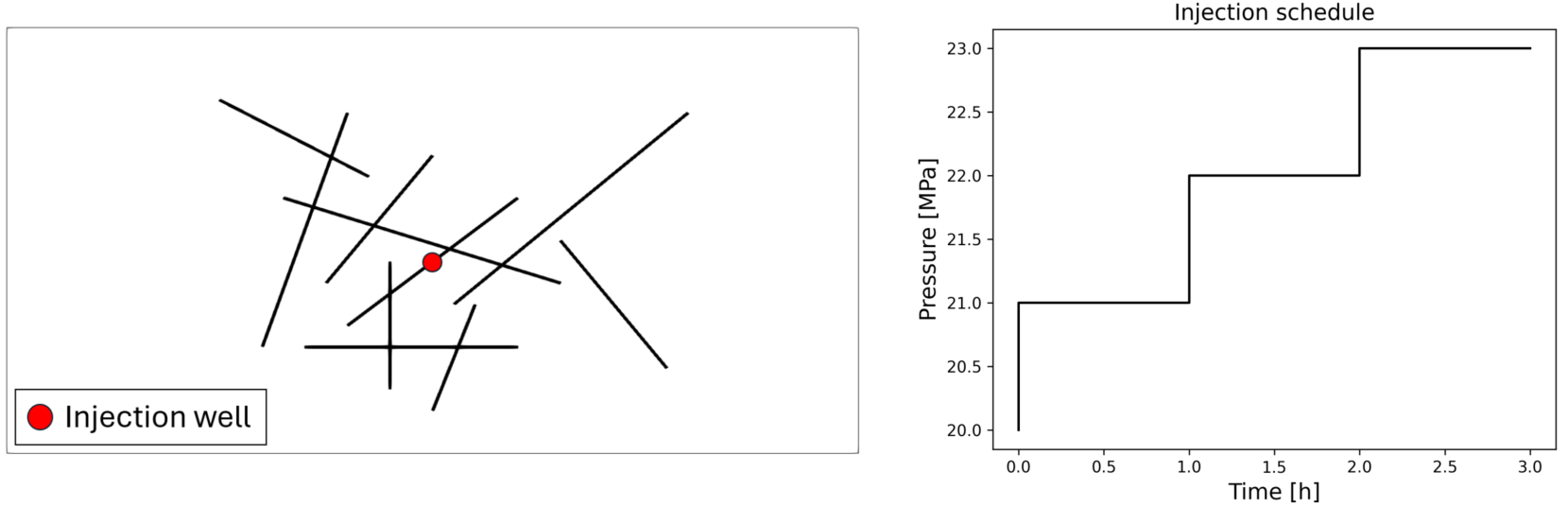}}}   
\caption{\change{Left: Fracture network used in the two-dimensional simulations, with the injection well marked by a red dot. Right: Injection schedule for the two-dimensional simulations. The wellhead pressure is momentarily increased by 1 MPa every hour, and otherwise kept constant. The same injection schedule will also be used for the three-dimensional simulations in section \ref{subsec:three_dim_ex}}}
    \label{fig:injection_schedule}
\end{figure}

\change{The performance of the nonlinear solvers, with the objectives of section \ref{sec:main_objectives} in mind, will first be explored in a series of two-dimensional numerical simulations.} A domain of size $2000 \times 1000$ m is considered, with the fracture network as shown in figure \ref{fig:injection_schedule}, and with a grid containing 33830 total cells (this includes matrix, fracture and interface cells). \change{Gravitational effects are neglected in these two-dimensional examples.} For the flow boundary conditions, we impose pressure values of 20MPa on all external boundaries. For the mechanical boundary conditions, we fix the bottom boundary (zero displacement), while compressive tractions are imposed on the remaining boundary sides. The traction values correspond to a stress tensor with non-zero values $\sigma_{xx}=-30$ MPa and $\sigma_{yy}=-50$ MPa. These pressures and anisotropic stress values are representative for crystalline rocks around 2 km depth. \change{The initial state is chosen to be the equilibrium state, with the pressure field being uniformly equal to 20 MPa.}

\change{The pressure in the centermost cell of one of the fractures (marked by a red dot in figure \ref{fig:injection_schedule}) will momentarily increase at certain points, and otherwise be kept at a constant pressure; this mimics a pressure-controlled injection well. The injection schedule is also shown in figure \ref{fig:injection_schedule}. Starting from the background pressure of 20 MPa, the overpressure in the well is increased by 1 MPa every hour, performing the first pressure increase at $t_{\text{start}}=0$, and running the simulation until $t_{\text{end}}=3$h. We will refer to the period between each pressure increase as the different \textit{phases} of the simulation. To accurately capture the dynamics immediately after each pressure increase, a small time step of one second is initially used. The adaptive time-stepping scheme described earlier is subsequently employed.}

\begin{figure}
[H]\centering
\includegraphics[width=\textwidth]{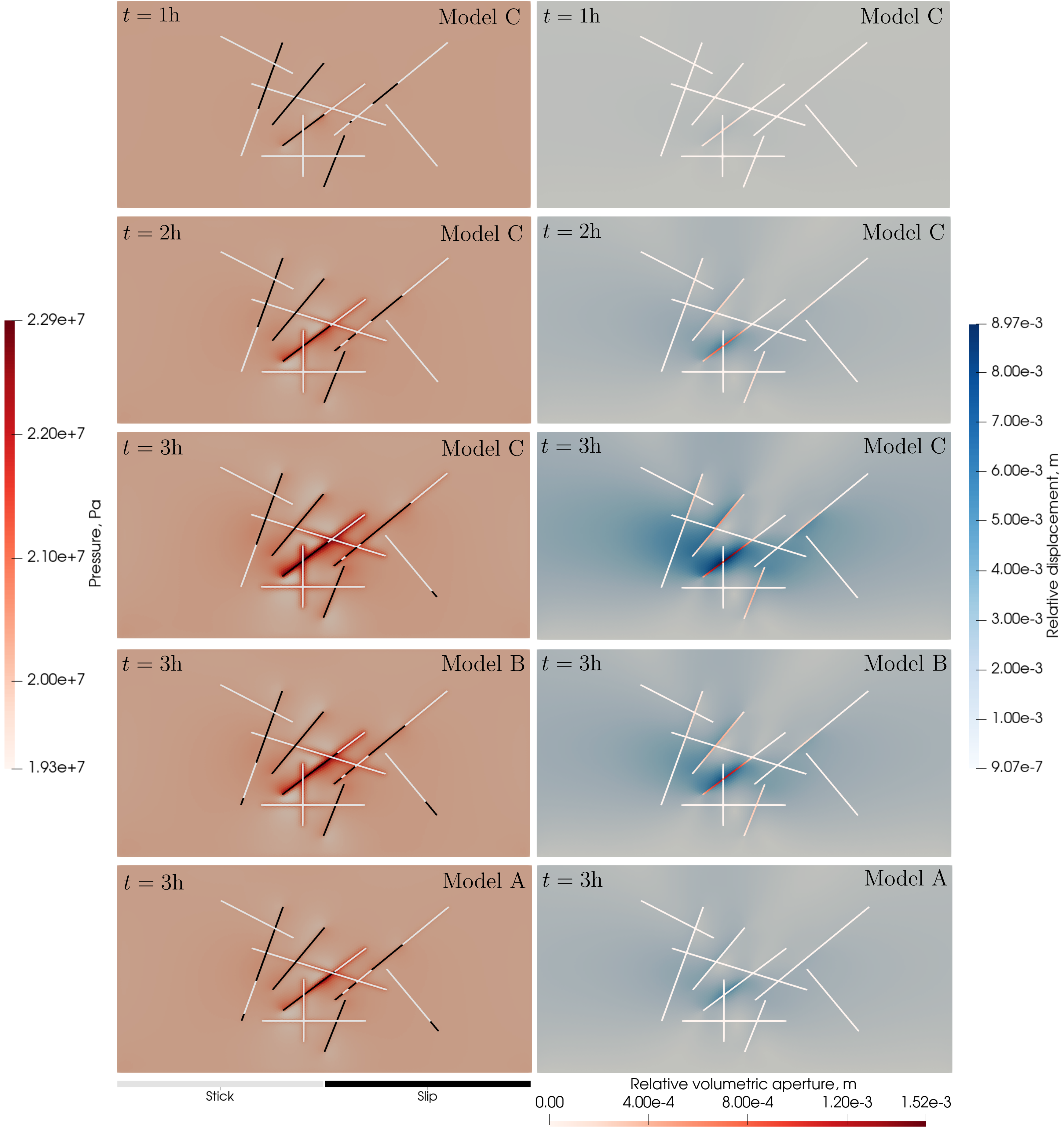}
\caption{\change{Visualization of the simulations of two-dimensional injection. Results are shown at the end of each simulation phase using model C, and at the end of the final simulation phase using models A and B. Left: Pressure profile and contact states. Here, slip is measured in a cumulative manner; a cell is regarded as being in a "slip" state if the fracture is closed and the norm of the tangential displacement jump is positive, using the state before the injection (at $t=0$) as the reference. There is no opening of fractures. Right: Changes in volumetric aperture and displacement magnitude, relative to the respective values at $t=0$}}
    \label{fig:pressure_states_2D}
\end{figure}

\change{We run simulations with the setup described above, using all three models of table \ref{tab:models_experiments}. Figure \ref{fig:pressure_states_2D} provides a visualization of the simulation at the end of each phase of the injection schedule using model C. For comparison, we also include in the figure the end of the final phase using models A and B. The plots show the pressure profile and contact states (left column), as well as the changes in volumetric aperture and displacement magnitudes relative to the initial state at $t=0$ (right column). We observe that each pressure increase in the well causes a larger portion of the fractures to slip. In the cases where models B and C are used, the slip events cause an increase in the volumetric aperture (this obviously does not happen in model A, where the volumetric aperture is constant). The increase in volumetric aperture in models B and C occurs due to the shear dilation of fractures causing an increase in the normal displacement jump. All models result in mostly the same fracture cells slipping, though with some minor differences. We also note that the pressures, displacement and volumetric apertures all gradually increase from models A through C.}

\begin{figure}[H]
    \centering
    \begin{subfigure}{0.49\textwidth} 
    \includegraphics[width=\textwidth]{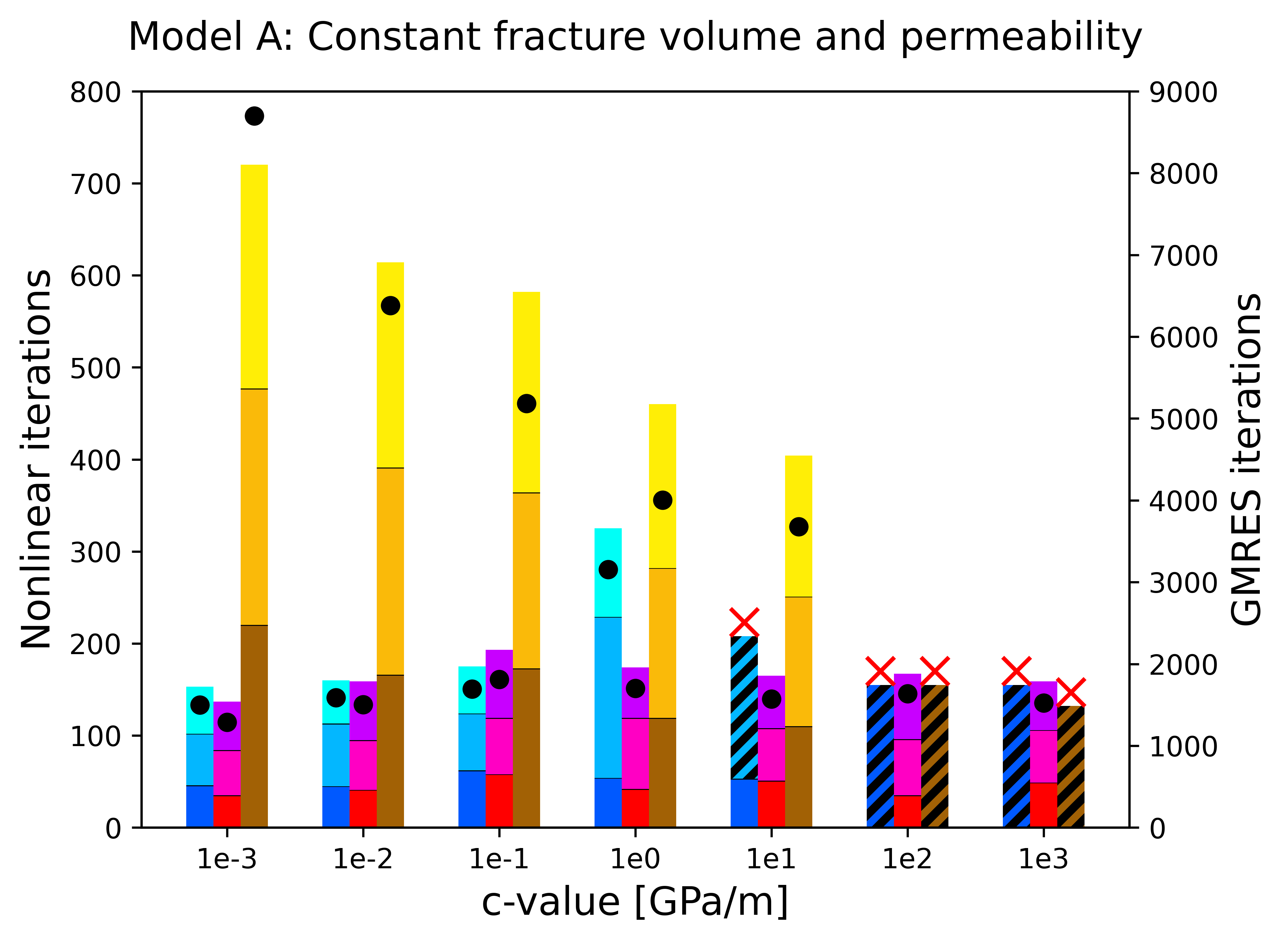}
    \caption{}
    \label{subfig:bar_no_aperture_2D}
    \end{subfigure}
    \begin{subfigure}{0.49\textwidth} 
    \includegraphics[width=\textwidth]{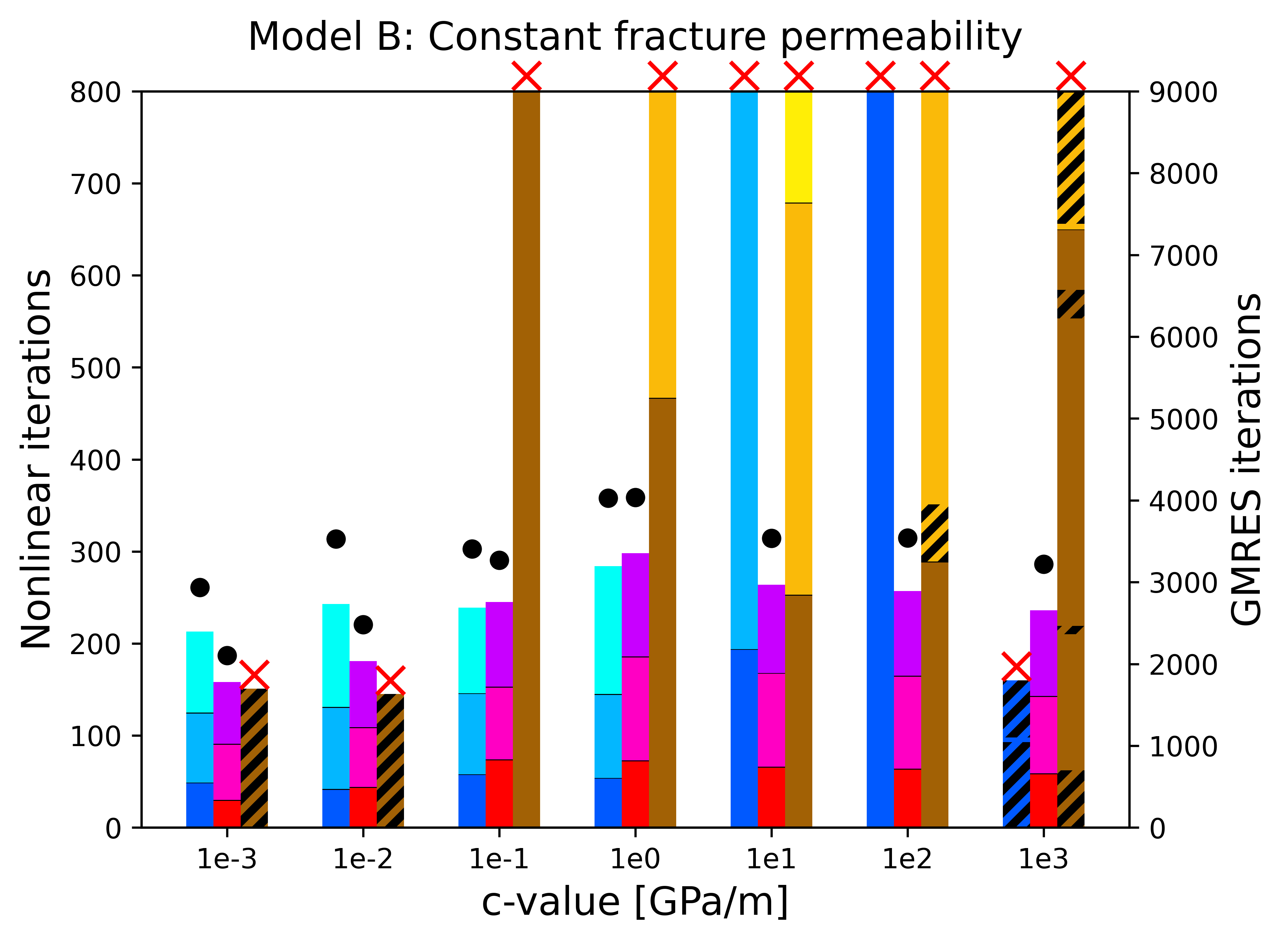}
    \caption{}
    \label{subfig:bar_no_cubic_2D}
    \end{subfigure}
    \begin{subfigure}{0.49\textwidth} 
    \includegraphics[width=\textwidth]{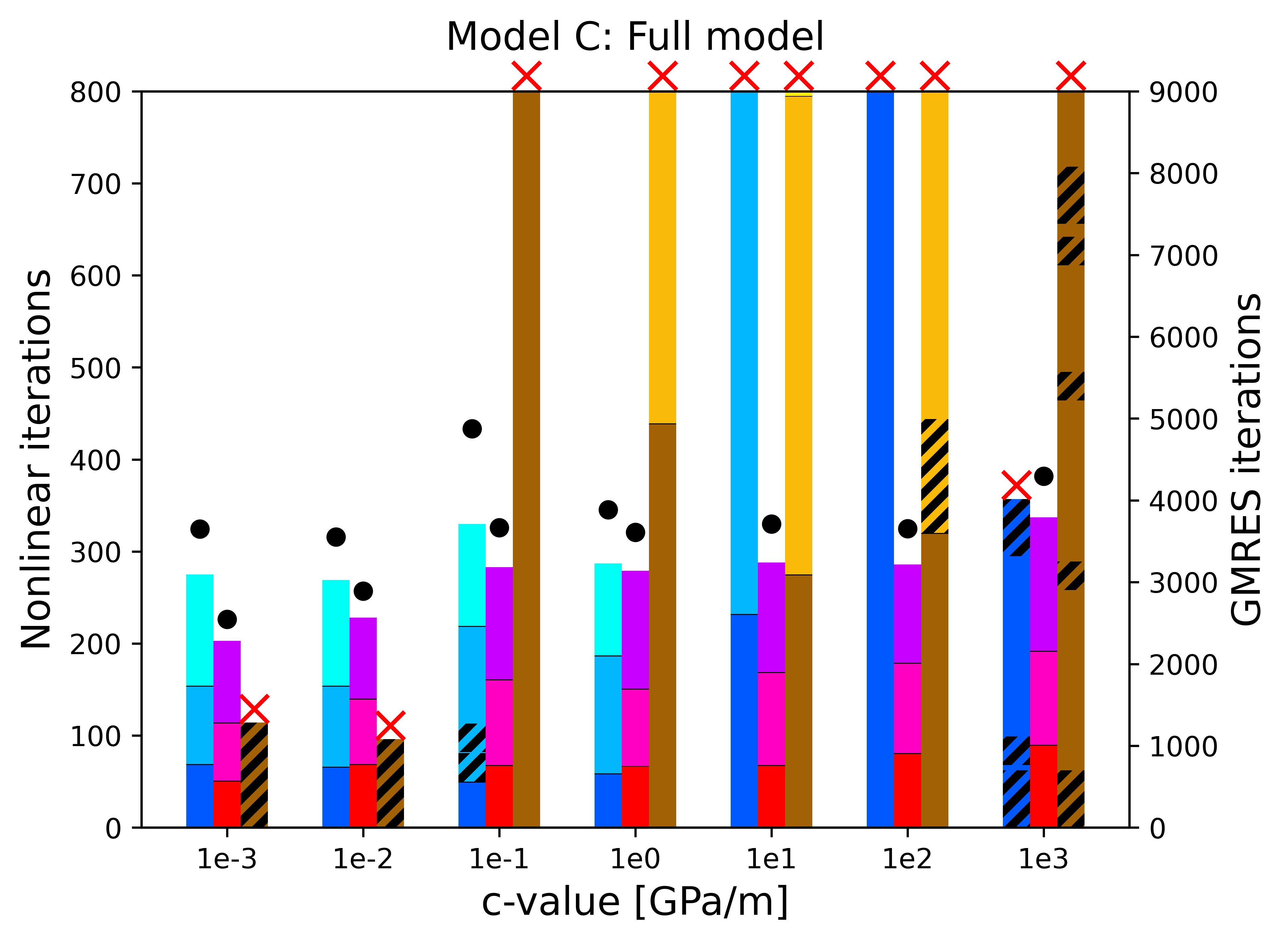}
    \caption{}
    \label{subfig:bar_full_model_2D}
    \end{subfigure}
    \begin{subfigure}{0.49\textwidth} 
    \includegraphics[width=\textwidth]{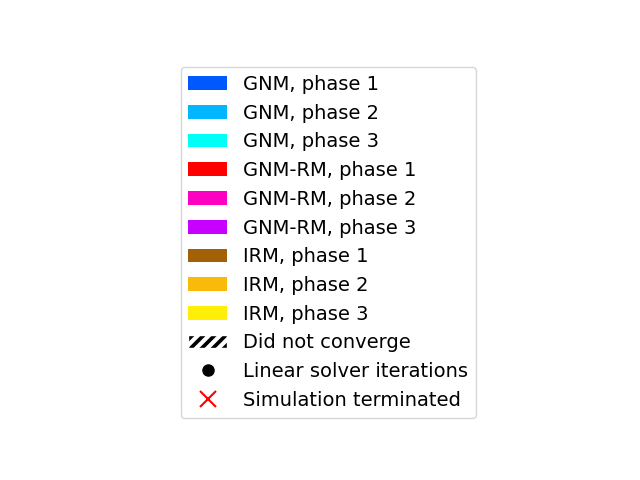}
    \caption{}
    \label{subfig:legend}
    \end{subfigure}
    \caption{\change{Bar plots showing the cumulative number of nonlinear iterations, summed over all time steps, for the example of section \ref{sec:two_dim_sim} with the different degrees of nonlinear couplings between flow and mechanics outlined in table \ref{tab:models_experiments}. The three different colors per bar correspond to the three different phases of the simulation, after each pressure increase in the injection well. Hatched black lines indicate steps where the nonlinear solver failed to converge, and red crosses indicate a simulation that was terminated. The simulation is terminated if one of three cases occur: The number of nonlinear iterations exceed 800, the nonlinear solver failed to converge for the smallest allowed time step, or after 6 unsuccessful recomputation attempts. Finally, the total number of iterations of the linear solver is marked by a black dot}}
    \label{fig:ex1_heatmaps}
\end{figure}

\change{In figure \ref{fig:ex1_heatmaps}, we report on the cumulative number of linear and nonlinear iterations, summed over all time steps, for the different nonlinear solvers and for different values of the augmentation parameter. As explained earlier, an update strategy is used on the augmentation parameter for IRM, and so the $c$-values in figure \ref{fig:ex1_heatmaps} refer to the starting values.}

\change{With the objectives of section \ref{sec:main_objectives} in mind, we start by discussing the impact of the different constitutive laws for the apertures. From figure \ref{fig:ex1_heatmaps}, we can clearly observe an increase in difficulty from model A through model C, with the number of iterations needed for convergence, as well as the number of non-convergent cases, generally increasing (note that the simulation is terminated if the number of nonlinear iterations exceeds 800; these cases may have converged if more iterations were used, but they are deemed too slow). The increase in difficulty is most pronounced between models A and B, with a sharp increase in the number of cases where the nonlinear solvers do not converge. We note in particular that IRM performs very poorly on models B and C, with no cases converging within 800 iterations, and several cases not even getting past the first injection phase. Overall, these results indicate that the nonlinear coupling between flow and mechanics, induced by the displacement-dependent aperture formulas, significantly adds to the difficulty of the problem, and that IRM is worse at handling this nonlinear coupling than the other two solvers.} 

\change{The performance of the linear solver is discussed next. We note from figure \ref{fig:ex1_heatmaps} a significant increase in the average number of GMRES iterations per nonlinear iteration from models A to B. This is consistent with the observation that the overall difficulty of the problem is significantly increased between these two models. Another interesting observation is that for IRM, the average number of GMRES iterations per nonlinear iteration generally decreases with increasing $c$-value, while it is roughly constant for GNM and GNM-RM ($c=10$ is an exception for IRM; here the number slightly increases again). Finally, it is worth noting that for models B and C, GNM-RM uses consistently less GMRES iterations per nonlinear iteration than GNM. Hence, for the cases considered here, the return map results in a better performance of the preconditioned linear solver.}

\change{Next, we discuss the impact of the augmentation parameter on the solver performances. First, we observe from figure \ref{subfig:bar_no_aperture_2D} the feature of IRM being slower to converge for lower $c$-values, which is present even with our relatively aggressive update strategy on the $c$-value. This feature, combined with the fact that IRM must solve more nonlinear systems per time step, causes IRM to be noticeably slower than GNM and GNM-RM. Second, we observe that both GNM and IRM perform poorly for the largest $c$-values in most cases, where they are either not converging at all, or using an excessive amount of nonlinear iterations. Additionally, IRM fails to converge also for the two smallest $c$-values for models B and C. GNM-RM, on the other hand, converges for all $c$-values, and is also quite consistent across the range of $c$-values, in regards to the number of nonlinear iterations. These results indicates that GNM-RM is less sensitive to the value of the augmentation parameter than the other two solvers, and in particular, that it handles higher values of the parameter much better.}

\begin{figure}[H]
    \centering
    \captionsetup[subfigure]{labelformat=empty}
    \begin{subfigure}{0.49\textwidth} 
    \includegraphics[width=\textwidth]{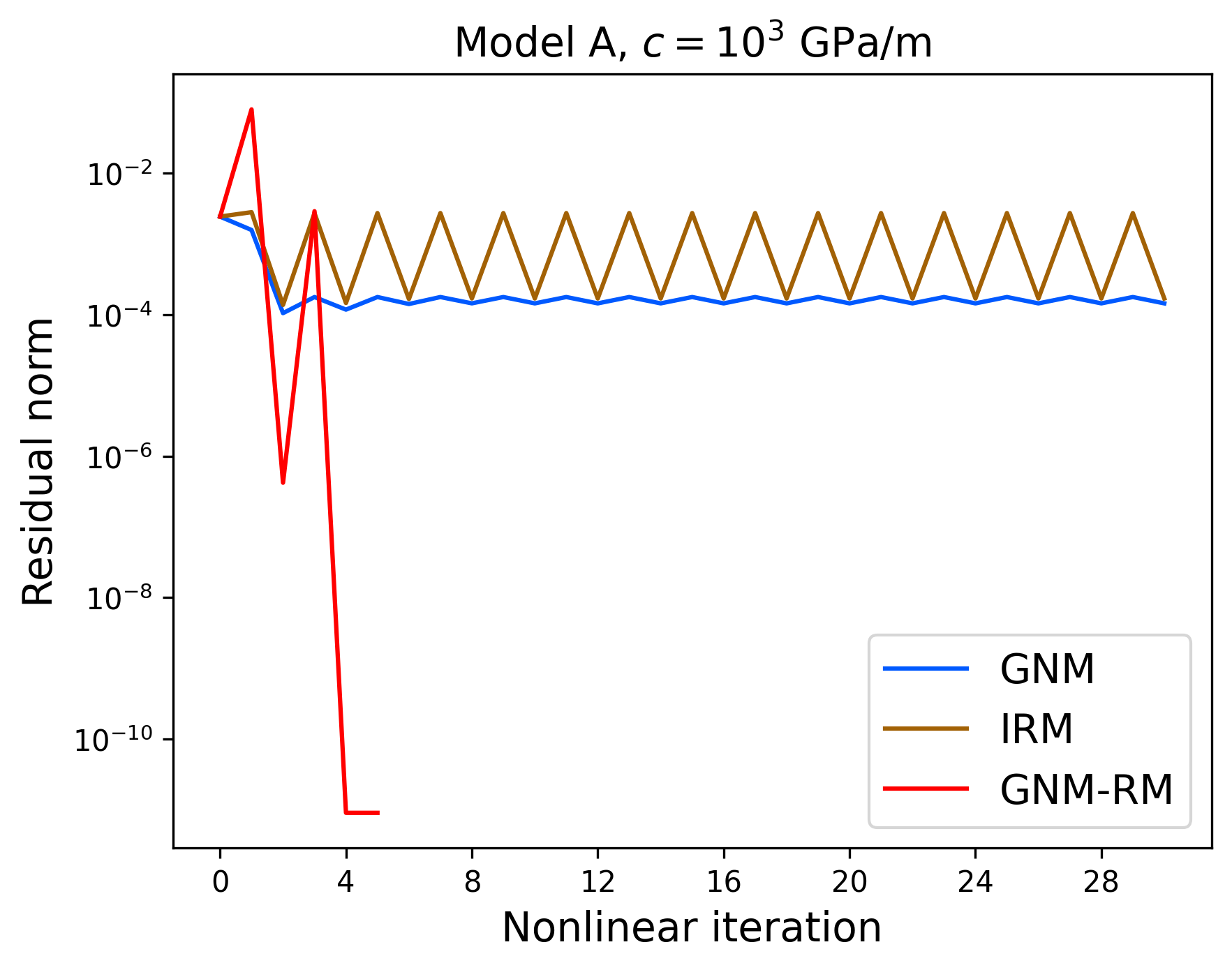}
    \caption{}
    \label{subfig:res_plots_cycle}
    \end{subfigure}
    \begin{subfigure}{0.49\textwidth} 
    \includegraphics[width=\textwidth]{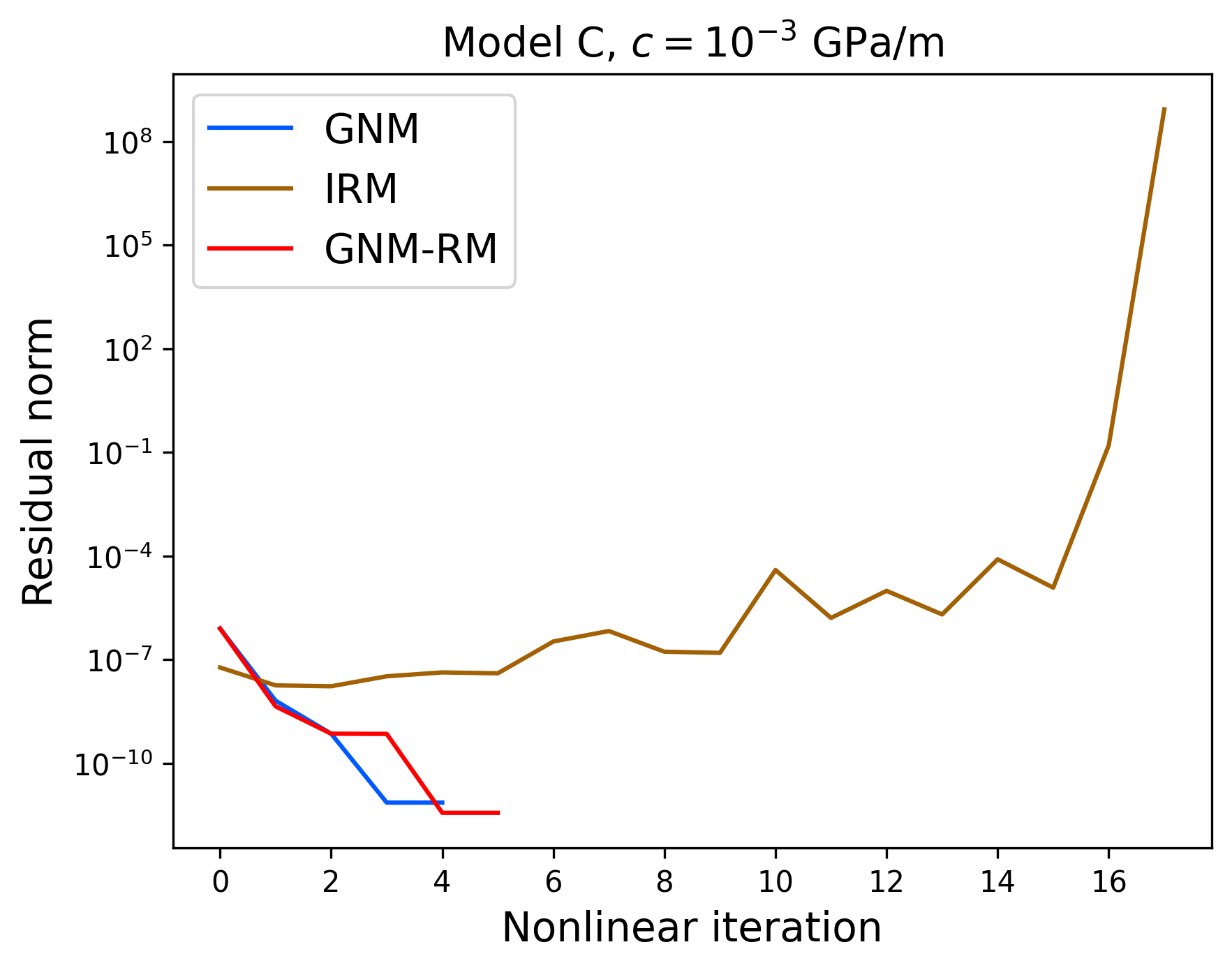}
    \caption{}
    \label{subfig:res_plots_div}
    \end{subfigure}
    \caption{\change{Residual norms at every nonlinear iteration for the simulation using model A and $c=10^3$ (left), and the simulation using model C and $c=10^{-3}$ (right). Only the first attempted time step is shown for all solvers}}
    \label{fig:ex1_residuals}
\end{figure}

\change{A closer inspection of the cases where GNM and IRM fail to converge (hatched lines in figure \ref{fig:ex1_heatmaps}) reveals a correlation between the mode of failure and the value of the augmentation parameter. For the higher values ($c>0$), the failure mode is nearly always cycling of the iterations, while for the lower values ($c<0$), the iterations either diverge to infinity or simply stagnate with no apparent pattern. Examples of both cases are shown in figure \ref{fig:ex1_residuals}, in which residual norms are plotted for the simulation using model A and a $c=10^3$, and the simulation using model C and $c=10^{-3}$. Different modes of failure are observed in the two cases. In the case with the higher $c$-value, GNM and IRM fails by cycling of the iterations, illustrated by repetitive patterns in the residual history. In the case with the lower $c$-value, IRM instead fails by the iterations diverging to infinity. GNM-RM is seen to converge quickly in both cases. An interesting conclusion to be drawn here is that, for these examples, performing a return map after each generalized Newton iteration is an effective way of avoiding cycling of the iterations.}

\change{In light of the entire comparison, we conclude that for these examples, IRM is the worst-performing solver, while GNM-RM is the best-performing one. IRM is both slower to converge than the other two solvers, as well as having more convergence issues for the models with the stronger nonlinear couplings between flow and mechanics. In terms of the number of nonlinear iterations, GNM and GNM-RM perform relatively similar. GNM-RM has the added benefit of being less sensitive to the value of the augmentation parameter than GNM.}

\change{The scalability of the solvers with respect to the grid resolution is also of interest. We performed an experiment testing the effect of the grid resolution on the number of nonlinear iterations needed for convergence. We have chosen not to include this experiment, however, as we were not able to draw any reliable conclusions from it. The general tendency was a weak grid dependency, with a modest increase in the number of iterations for finer grids in most cases. However, this pattern was also broken on several occasions, as the grid refinement would sometimes result in a fewer number of iterations.}

\subsection{Three-dimensional simulations of fluid injection} \label{subsec:three_dim_ex}

\begin{figure}
[H]\centering
\scalebox{0.6}{\includegraphics[width=\textwidth]{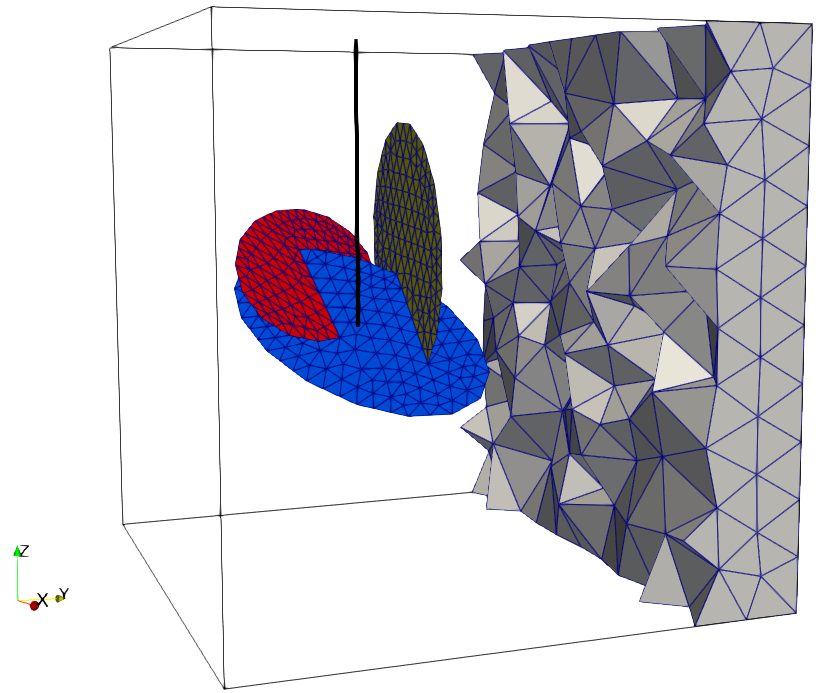}}
\caption{\change{Fracture network geometry used in the three-dimensional simulations, and the position of the injection well, marked by the vertical black line. Also shown are the two-dimensional grid cells, as well as parts of the three-dimensional grid. The grid contains a total of 16405 cells, including matrix, fracture and interface cells. The average size of the grid cells is seen to be larger near the boundary than on the fractures, indicating grid refinement near the fractures.
The different colors on the fractures are used to visually distinguish them, and do not refer to any physical quantity}}
\label{fig:setup_ex_3}
\end{figure}

\change{Next, we investigate the performance of the nonlinear solvers on a set of three-dimensional problems, and investigate to which degree the results from the two-dimensional case carry over to three dimensions. As before, our aim is} to simulate hydraulic stimulation of critically stressed fractures by an injection well. The domain is the box $(0 \text{m}, 2000 \text{m}) \times (0 \text{m}, 2000 \text{m}) \times (-3000 \text{m}, -1000 \text{m})$, and contains a network of three intersecting, elliptic-shaped fractures, as shown in figure \ref{fig:setup_ex_3}. For the hydraulic boundary conditions, we impose hydrostatic pressure conditions on all external boundaries, given by $p=\rho^f_{\mathrm{ref}} g z$, with $z$ being the depth and $g$ the gravitational acceleration. For the mechanical boundary conditions, we fix the bottom boundary (zero displacement), while for the remaining boundaries, we impose compressive tractions, with a similar anisotropy as used in section \ref{sec:two_dim_sim}. Namely, the corresponding non-zero stress tensor values are:
\begin{equation}
    \sigma_{xx}=-0.6\rho^s g z \ \ , \ \ \sigma_{yy}=-0.6\rho^s g z \ \ , \ \ \sigma_{zz}=-\rho^s g z.
\end{equation}
The rest of the simulation setup is similar to that of section \ref{sec:two_dim_sim}. \change{The pressure in the centermost cell of one of the fractures is momentarily increased at certain points, mimicking an injection well (illustrated by the vertical black line in figure \ref{fig:setup_ex_3}). The total simulation time is three hours as in section \ref{sec:two_dim_sim}, with the pressure profile of the injection cell again following figure \ref{fig:injection_schedule}. The adaptive time-stepping scheme outlined in section \ref{sec:technical_details} is employed. The initial state is the equilibrium state with a hydrostatic pressure field. Again we will test the effects of different degrees of nonlinear couplings between flow and mechanics, by using the three different models of table \ref{tab:models_experiments}.}

\begin{figure}
[H]\centering
\scalebox{1}{\includegraphics[width=\textwidth]{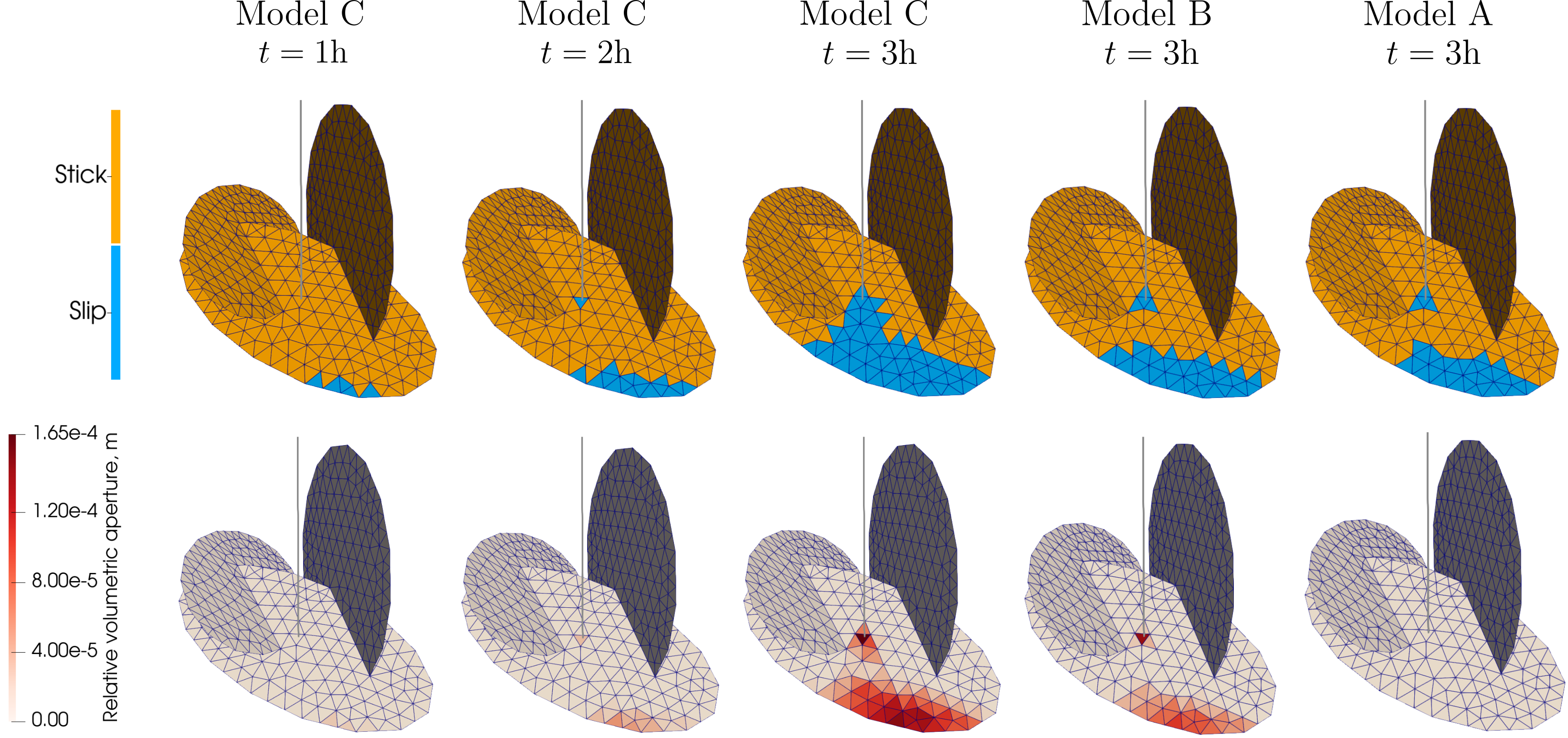}}
\caption{\change{Visualization of the simulations of three-dimensional injection. Results are shown at the end of each simulation phase using model C, and at the end of the final simulation phase using models A and B.
Top row: Contact states. Here, slip is measured in a cumulative manner; a cell is regarded as being in a "slip" state if the fracture is closed and the norm of the tangential displacement jump is positive, using the state before the injection (at $t=0$) as the reference. There is no opening of fractures. Bottom row: Changes in volumetric aperture, relative to the value at $t=0$}}
    \label{fig:states_aperture_3D}
\end{figure}

\change{The resulting contact states and differences in volumetric aperture is shown at various points in the simulations in figure \ref{fig:states_aperture_3D}. We note that slip around the injection well does not occur until the second phase, when a sufficiently high pressure is reached. Otherwise, we observe qualitatively similar results to the two-dimensional case. The slip events result in an increase of the volumetric aperture (except in the constant-aperture case of model A), due to shear dilation, and the magnitude of the volumetric aperture difference increases from models B to C.}

\begin{figure}[H]
    \centering
    \begin{subfigure}{0.49\textwidth} 
    \includegraphics[width=\textwidth]{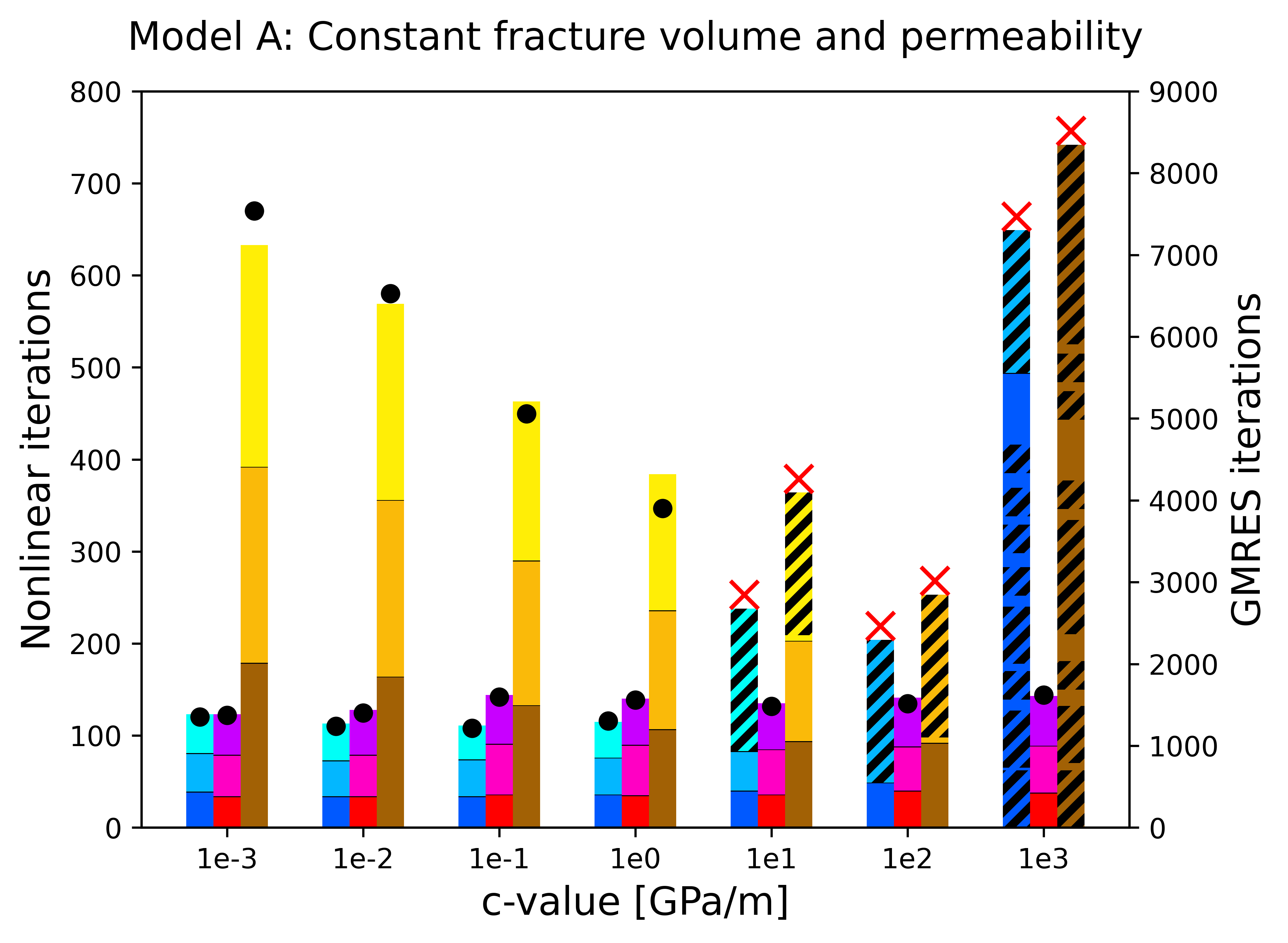}
    \caption{}
    \label{subfig:bar_no_aperture_3D}
    \end{subfigure}
    \begin{subfigure}{0.49\textwidth} 
    \includegraphics[width=\textwidth]{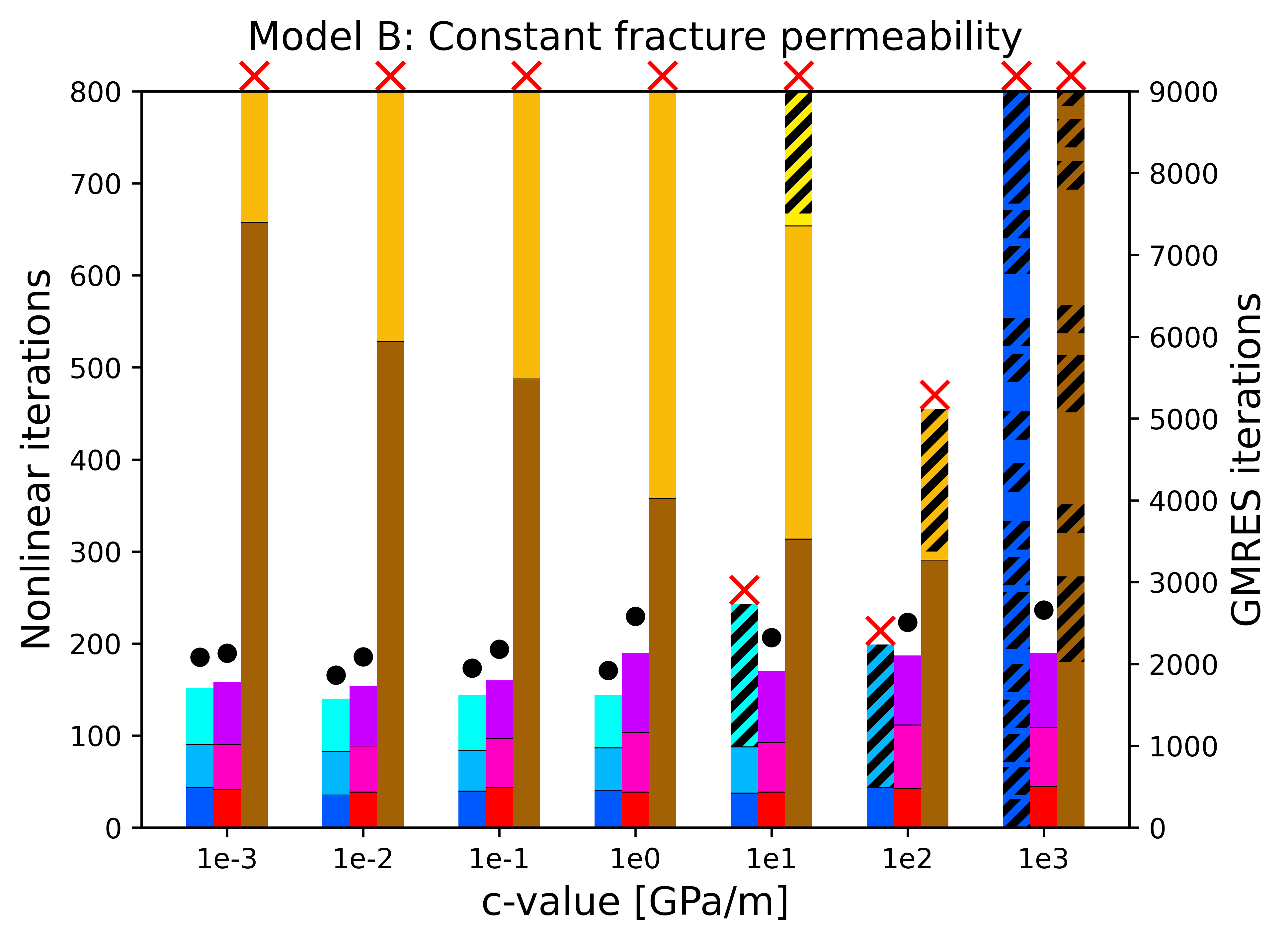}
    \caption{}
    \label{subfig:bar_no_cubic_3D}
    \end{subfigure}
    \begin{subfigure}{0.49\textwidth} 
    \includegraphics[width=\textwidth]{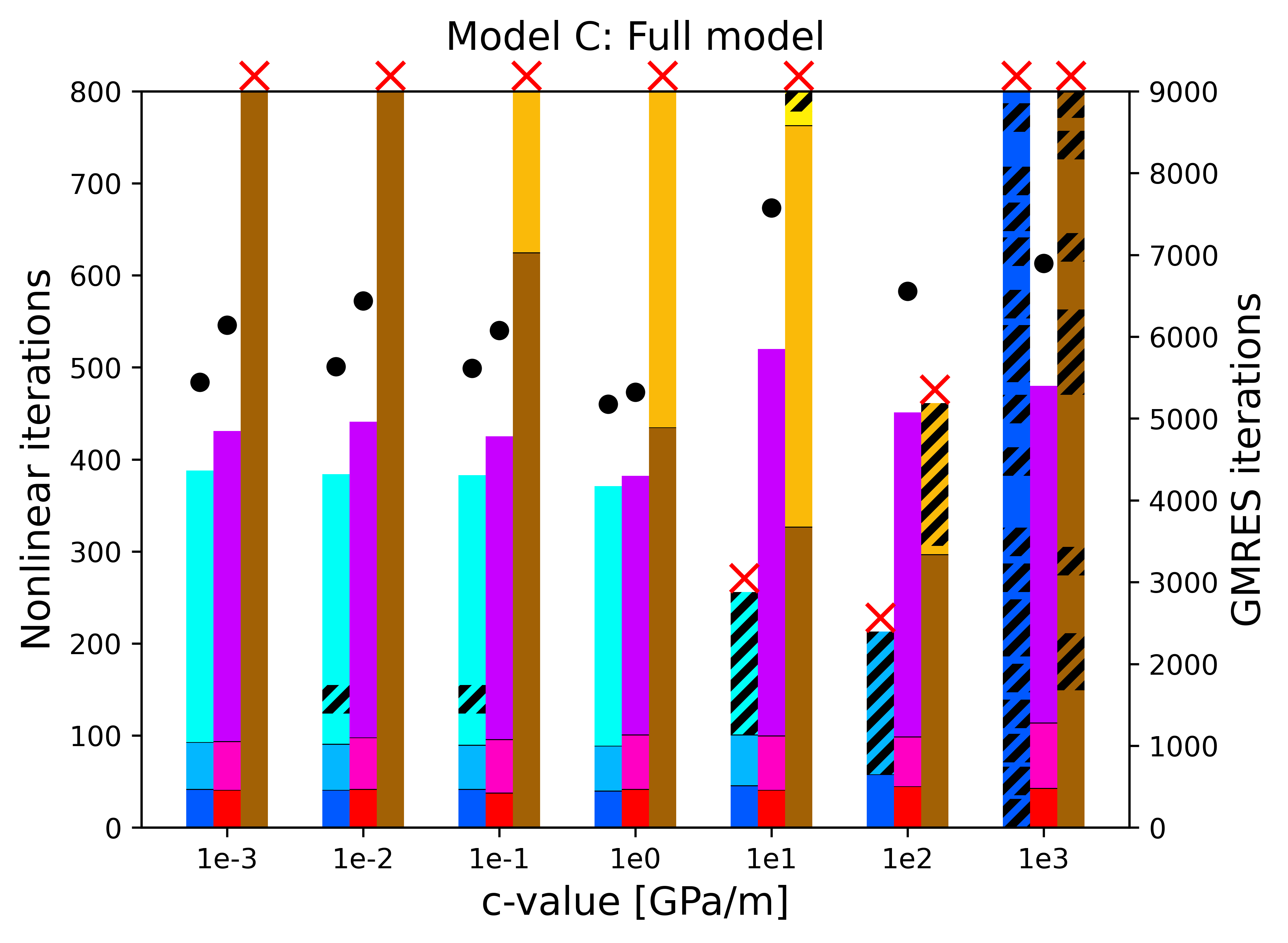}
    \caption{}
    \label{subfig:bar_full_model_3D}
    \end{subfigure}
    \begin{subfigure}{0.49\textwidth} 
    \includegraphics[width=\textwidth]{Fig4d.png}
    \caption{}
    \label{subfig:legend_again}
    \end{subfigure}
    \caption{\change{Bar plots showing the cumulative number of nonlinear iterations, summed over all time steps, for the example of section \ref{subsec:three_dim_ex} with the different degrees of nonlinear couplings between flow and mechanics outlined in table \ref{tab:models_experiments}. The three different colors per bar correspond to the three different phases of the simulation, after each pressure increase in the injection well. Hatched black lines indicate steps where the nonlinear solver failed to converge, and red crosses indicate a simulation that was terminated. The simulation is terminated if one of three cases occur: The number of nonlinear iterations exceed 800, the nonlinear solver failed to converge for the smallest allowed time step, or after 6 unsuccessful recomputation attempts. Finally, the total number of iterations of the linear solver is marked by a black dot}}
    \label{fig:bar_charts_3D}
\end{figure}

\change{The performance of the nonlinear solvers for different values of the augmentation parameter is shown in figure \ref{fig:bar_charts_3D}. We observe many of the same patterns as in the two-dimensional simulations. An increase in difficulty occurs from models A through C, as the nonlinear coupling strength is increased. The performance of IRM deteriorates for models B and C, with no cases converging within 800 iterations. In regards to the augmentation parameter, GNM and IRM both fail to converge for the largest three values, while GNM-RM converges in all cases. These results are all very similar to the results of section \ref{sec:two_dim_sim}.}

\change{The behavior of the linear solver also mostly matches that of the two-dimensional simulations. The average number of GMRES iterations per nonlinear iteration slightly decreases with increasing augmentation parameter for IRM, while it is close to constant with respect to the augmentation parameter for GNM and GNM-RM. The same behavior was seen in the two-dimensional simulations of figure \ref{fig:ex1_heatmaps}. A slight difference from the two-dimensional simulations is that for the cases considered here, the return map of GNM-RM does not appear to improve the performance of the linear solver. In nearly all cases of figure \ref{fig:ex1_heatmaps}, GNM-RM used less GMRES iterations per nonlinear iteration than GNM, while in figure \ref{fig:bar_charts_3D}, this number is about equal for GNM and GNM-RM.}

\change{Another notable difference from the two-dimensional simulations is that the fully coupled model is significantly more challenging than the other two models. It is apparent from figure \ref{subfig:bar_full_model_3D} that for model C, a considerable spike in difficulty occurs during the third phase, with a large increase in the number of nonlinear iterations. Additionally, GNM fails to converge even for some of the lower $c$-values during this phase (although it is successful after the time step is lowered).}

\subsection{More challenging three-dimensional simulations} \label{sec:three_dim_difficult}

\begin{figure}[H]\centering
\scalebox{0.6}{\includegraphics[width=\textwidth]{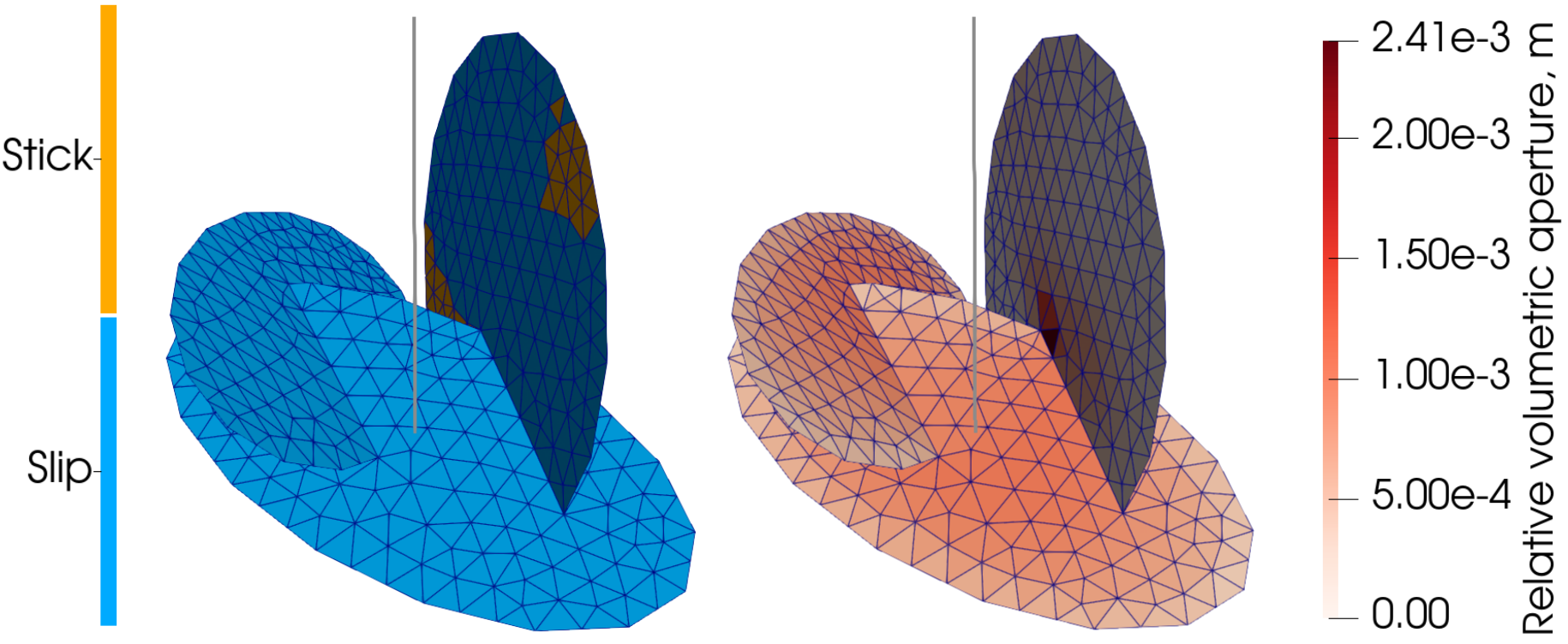}}
\caption{\change{Visualization of the simulation of section \ref{sec:three_dim_difficult}, after one hour of simulation time. Left: Fracture contact states. Here, slip is measured in a cumulative manner; a cell is regarded as being in a "slip" state if the fracture is closed and the norm of the tangential displacement jump is positive, using the state before the injection (at $t=0$) as the reference. There is no opening of fractures. Right: Changes in volumetric aperture, relative to the value at $t=0$}}
\label{fig:three_dim_states_hard}
\end{figure}

\begin{figure}[H]
    \centering
    \begin{subfigure}{0.49\textwidth}   \includegraphics[width=\textwidth]{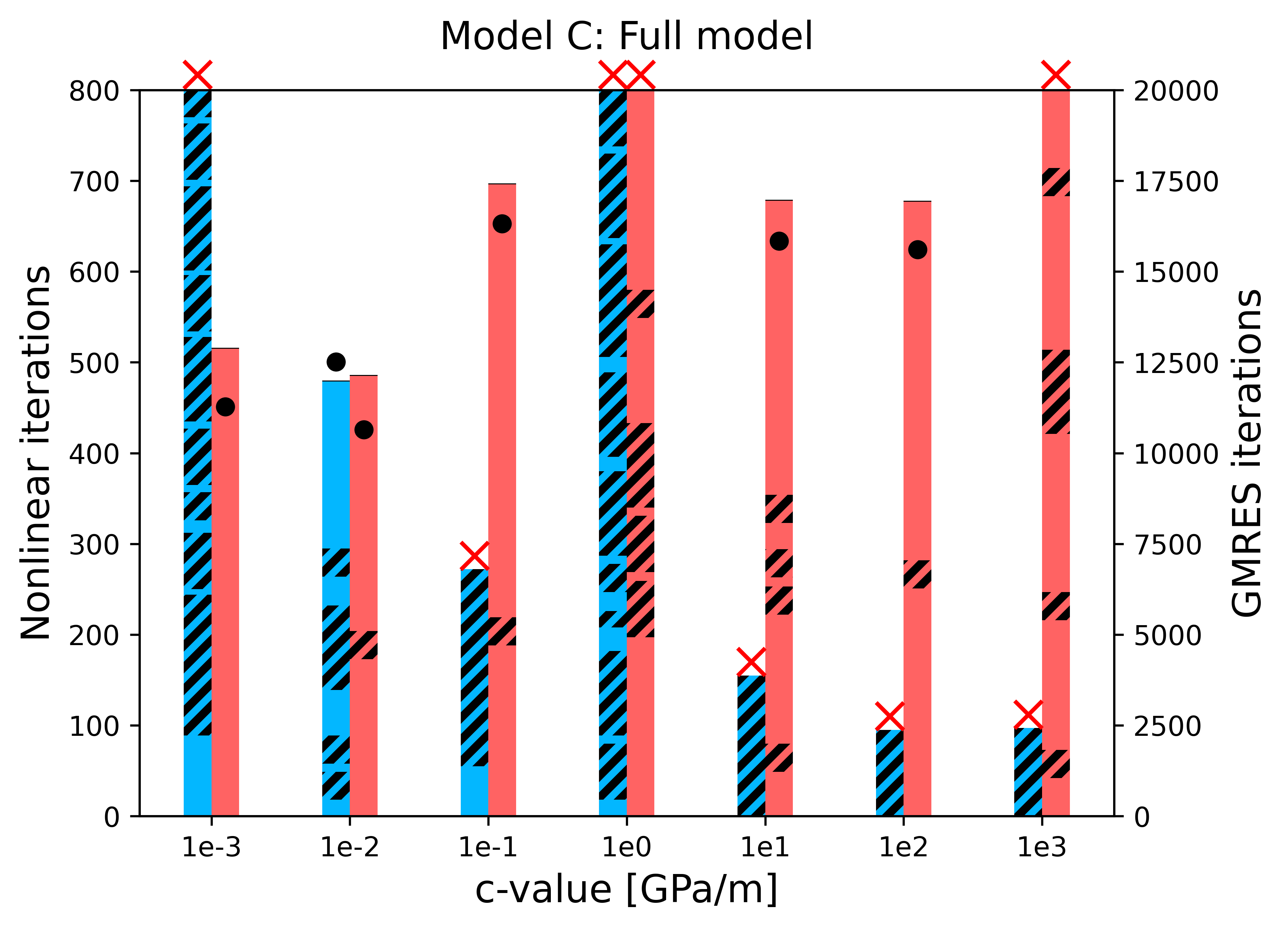}  \caption{}\label{subfig:bar_full_model_difficult}
    \end{subfigure}
    \begin{subfigure}{0.49\textwidth} 
    \includegraphics[width=\textwidth]{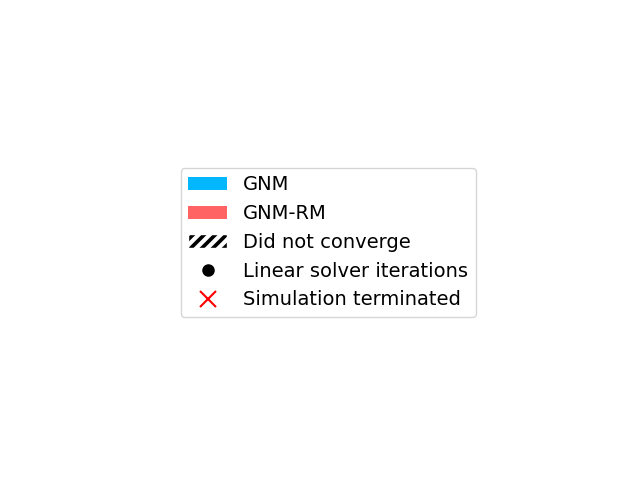}
    \caption{}
    \label{subfig:legend_difficult}
    \end{subfigure}
    \caption{\change{Bar plots showing the cumulative number of nonlinear iterations, summed over all time steps, for the example of section \ref{sec:three_dim_difficult}. Hatched black lines indicate steps where the nonlinear solver failed to converge, and red crosses indicate a simulation that was terminated. The simulation is terminated if one of three cases occur: The number of nonlinear iterations exceed 800, the nonlinear solver failed to converge for the smallest allowed time step, or after 6 unsuccessful recomputation attempts. Finally, the total number of iterations of the linear solver is marked by a black dot.}}
    \label{fig:bar_chart_difficult}
\end{figure}

\change{To further explore the challenges of fully coupled three-dimensional simulations and test the limits of GNM and GNM-RM, a final set of simulations has been run using model C. We discard IRM in this final set of simulations, as it has been consistently outclassed by the other two solvers, and instead focus on the comparison between GNM and GNM-RM. The simulation setup is nearly identical to the previous one; the only difference is a reduction in the compressive tractions on the boundary faces aligned with the $y$-direction. The stress field on the boundary is now as follows:}
\begin{equation} \label{eq:new_stress_field}
    \sigma_{xx}=-0.6\rho^s g z \ \ , \ \ \sigma_{yy}=-0.4\rho^s g z \ \ , \ \ \sigma_{zz}=-\rho^s g z.
\end{equation}
\change{This reduction in the compressive stress field causes a much larger portion of the fractures to slip, as shown in figure \ref{fig:three_dim_states_hard}, creating more complex nonlinear dynamics that challenge the solvers. Due to the added difficulty of this setup, we will run only the first phase of the injection schedule of figure \ref{fig:injection_schedule}, i.e. we increase the pressure in the injection well by 1MPa at $t=0$, and run the simulation for one hour. The results are reported in figure \ref{fig:bar_chart_difficult}.}

\change{We observe a significant increase in difficulty for the solvers. The total number of nonlinear iterations is now comparable to the previous three-dimensional simulation with model C (figure \ref{subfig:bar_full_model_3D}) for only one third of the simulation time. Moreover, in all cases the nonlinear solvers failed to converge at some point, resulting in a recomputation with a smaller time step. The number of GMRES iterations per nonlinear iteration is also much higher for this example than the previous ones (note the change of scale on the axis for the linear iterations, compared to the previous bar plots). Despite the increased difficulty, GNM-RM remains relatively robust, managing to converge within 800 iterations in all but two cases.}

\change{Overall, the results of the three-dimensional simulations of sections \ref{subsec:three_dim_ex} and \ref{sec:three_dim_difficult} are relatively similar to the two-dimensional simulations of section \ref{sec:two_dim_sim}. The increase in difficulty from models A to C is once again clearly seen. IRM is the worst-performing solver also in three dimensions, and deteriorates for models B and C, while GNM-RM is the best-performing one, being the most robust solver with respect to the augmentation parameter. However, we also observe some more unpredictable behavior of the nonlinear solvers for the most difficult case of section \ref{sec:three_dim_difficult}. This is unsurprising, given the high complexity of the fully coupled, three-dimensional model.}

\section{Conclusion} \label{sec:conclusion}

\change{We have investigated the performance of three different nonlinear solvers, all based on the augmented Lagrangian formulation of frictional contact mechanics, for solving advanced problems of coupled flow and deformation in fractured porous media, including fracture contact mechanics. This included two well-known solvers from the contact mechanics literature, namely the generalized Newton method (GNM) using complementarity functions, and the implicit return map method (IRM), the latter of which is equivalent to an implicit Uzawa method. In addition, we proposed a third nonlinear solver, which combined features of both GNM and IRM. This solver, named the ``generalized Newton method with a return map'' (GNM-RM), incorporates the return map into the generalized Newton method as a postprocessing step after every iteration. The purpose of the return map in GNM-RM is to stabilize the generalized Newton method, which is in contrast to classical return map methods like IRM, where its purpose is to ensure the feasibility of the contact tractions.}

\change{The performances of the three solvers were assessed through a series of numerical experiments, both two- and three-dimensional, which were designed to mimic hydraulic stimulation of a fractured geothermal reservoir. In the experiments, the degree of nonlinear coupling between flow and mechanics was varied, and different values of the augmentation parameter were tested. The experiments showed that increasing the nonlinear coupling between flow and mechanics significantly adds to the difficulty of the problem, with a higher frequency of non-convergent behavior of the nonlinear solvers, as well as a general increase in the number of nonlinear iterations needed for convergence. GNM and GNM-RM still managed to converge in many cases for the more strongly coupled problems, but IRM performed much worse, barely converging at all except for the cases with very weak couplings. In regards to the augmentation parameter, GNM-RM was the most robust solver, converging across the largest range of values of this parameter. In particular, it handled higher values of the parameter much better than the other two solvers. For the higher values of the augmentation parameter, GNM and IRM would frequently fail to converge, with cycling of the iterations being the central failure more, while GNM-RM would avoid the cycle and converge in nearly all cases. This makes GNM-RM a viable alternative to GNM for the type of problems considered herein.}

\section*{Acknowledgements}

\change{This project has received funding from the European Research Council (ERC) under the
European Union’s Horizon 2020 research and innovation program (grant agreement No.
101002507).}

\bibliography{bibliography}

\end{document}